%% file: paper.tex
\documentclass[11pt, reqno]{amsart}

\usepackage {amsmath, amssymb, amscd, mathrsfs, url,verbatim,enumerate}
\usepackage[text={6.5in,9in},centering,letterpaper,dvips]{geometry}
\usepackage{multirow,latexsym,graphicx,epstopdf,comment}
\usepackage[all]{xy}

\usepackage[dvipsnames,usenames]{color}
\usepackage[colorlinks=true, urlcolor=NavyBlue, linkcolor=NavyBlue, citecolor=NavyBlue]{hyperref}
\usepackage{graphicx}
\usepackage{epsfig}
\usepackage[latin1]{inputenc}
\usepackage{amsmath}
\usepackage{amsfonts}
\usepackage{amssymb}
\usepackage{amsthm}
\usepackage{amscd}
\usepackage{verbatim}
\usepackage{subfigure}
\usepackage{pinlabel}
\usepackage{stmaryrd}
\usepackage{enumerate, enumitem}

\usepackage{tikz}
\usetikzlibrary{arrows}
\usetikzlibrary{decorations.pathreplacing}
\usepackage{verbatim}
	
\newtheorem{theorem}{Theorem}[section]

\newtheorem{lemma}[theorem]{Lemma}
\newtheorem{proposition}[theorem]{Proposition}

\newtheorem {fact}[theorem]{Fact}

\theoremstyle{definition}
\newtheorem{definition}[theorem]{Definition}

\theoremstyle{remark}
\newtheorem{remark}[theorem]{Remark}

\def\gr{\mathrm{gr}}

\def\Z{\mathbb{Z}}

\def\F{\mathbb{F}}

\def\min{\textup{min}}

\def\Tplus{\mathcal{T}^+}
\newcommand{\HFred}{HF^{\textup{red}}}

\DeclareMathOperator{\rank}{rank}
\newcommand{\Spinc}{$\mathrm{Spin}^c\;$}

\DeclareMathOperator{\Sym}{Sym}
\DeclareMathOperator{\suc}{suc}
\DeclareMathOperator{\pre}{pre}

\def\red{\textup{red}}
\def\deg{\textup{deg}}
\def\Im{\textup{Im}}

\specialcomment{jen}{%
  \begingroup\color{green}  Jen\; $\blacktriangleright$ \;}{%
	\endgroup}
\specialcomment{cagri}{%
  \begingroup\color{blue}  Cagri \; $\blacktriangleright$ \;}{%
  \endgroup}
\specialcomment{tye}{%
  \begingroup\color{red} Tye \; $\blacktriangleright$ \;}{%
  \endgroup}

\subjclass[2013]{}
\author[Jennifer Hom]{Jennifer Hom}
\thanks{The first author was partially supported by NSF grant DMS-1307879.}
\address {Department of Mathematics, Columbia University, New York, NY 10027}
\email{hom@math.columbia.edu}

\author[\c{C}a\u{g}r{\i} Karakurt]{\c{C}a\u{g}r{\i} Karakurt}
\thanks{The second author was partially supported by NSF grant DMS-1065178.}
\address{Department of Mathematics, Bo{\u{g}}azi{\c{c}}i University, Bebek, {\.{I}}stanbul, TR-34342, Turkey, and, The University of Texas, Austin, TX 78712}
\email{karakurt@math.utexas.edu}

\author[Tye Lidman]{Tye Lidman}
\thanks{The third author was partially supported by NSF grant DMS-0636643.}
\address {Department of Mathematics, The University of Texas, Austin, TX 78712}
\email{tlid@math.utexas.edu}

\numberwithin{equation}{section}

\title{Surgery obstructions and Heegaard Floer homology}
\date{}

\begin{document}
\maketitle

\begin{abstract}
Using Taubes' periodic ends theorem, Auckly gave examples of toroidal and hyperbolic irreducible integer homology spheres which are not surgery on a knot in the three-sphere.  We give an obstruction to a homology sphere being surgery on a knot coming from Heegaard Floer homology.  This is used to construct infinitely many small Seifert fibered examples.  
\end{abstract}

\input{intro}
\input{mappingcone}
\input{plumbings}

\input{graded}
\input{delta}

\input{semigroups}
\input{universalspheres}

\bibliography{References}
\bibliographystyle{amsalpha}

\end{document}

%% file: intro.tex
\section{Introduction}\label{sec:intro}

\subsection{Background}
A classical theorem due to Lickorish and Wallace states that every closed oriented three-manifold can be expressed as surgery on a link in $S^3$ \cite{Lickorish,Wallace}.  Therefore, a natural question is which three-manifolds have the simplest surgery presentations.  More specifically, we ask which three-manifolds can be represented by surgery on a knot.  There are a number of obstructions that one can apply (with a range of effectiveness) to attempt to answer this question.  Since $H_1(S^3_{p/q}(K);\mathbb{Z}) \cong \mathbb{Z}/p$, one finds immediate homological obstructions to a manifold being obtained by surgery on a knot (e.g., $\mathbb{T}^3$ and $\mathbb{R}P^3 \# \mathbb{R}P^3$ are not surgery on a knot).  A more delicate obstruction is the weight of the fundamental group; a three-manifold obtained by surgery on a knot in $S^3$ has weight one fundamental group (normally generated by a single element).  Observe that weight one groups necessarily have cyclic abelianization.  Hence, this obstruction extends the aforementioned homological obstruction.  

A more topological obstruction can be found in the prime decomposition of the three-manifold.  A theorem of Gordon and Luecke \cite{GL} shows that if surgery on a non-trivial knot in $S^3$ yields a reducible manifold, one of the summands is necessarily a non-trivial lens space.  In particular, a reducible integer homology sphere can never be surgery on a knot.  Boyer and Lines \cite{BL} give an infinite family of prime Seifert fibered manifolds with weight one fundamental group which are not surgery on a knot.  Their proof requires two obstructions: the first comes from an extension of the Casson invariant to homology lens spaces and the second comes from the linking form.  In particular, having non-trivial first homology was necessary in their proof.                

Utilizing Taubes' periodic ends theorem \cite{T}, Auckly constructed examples of irreducible integer homology spheres which are not surgery on a knot in $S^3$, answering \cite[Problem 3.6(C)]{K}.  The first example was toroidal and shown to be homology cobordant to $\Sigma(2,3,5) \# - \Sigma(2,3,5)$ (or equivalently, $S^3$).  This was then extended to give a hyperbolic example \cite{Auckly}.  However, as far as the authors know, it is unknown whether Auckly's examples have weight one fundamental group.  This is related to the question of Wiegold (see, for instance, \cite[Problem 5.52]{Groups}): does every finitely-presented perfect group have weight one?  

In \cite{S}, it is asked if there are Seifert fibered homology spheres which are not surgery on a knot.  Note that every Seifert homology sphere is irreducible, so none of these are ruled out by the Gordon-Luecke criterion.  We answer this question affirmatively.

\subsection{Main results}

\begin{theorem}\label{thm:main}
For $p$ an even integer at least 8, let $Y_p$ denote the Seifert fibered integer homology sphere $\Sigma(p,2p-1,2p+1)$.  The manifolds $Y_p$ satisfy:   
\begin{enumerate}[label=(\roman{*}), ref=\roman{*}]
\item \label{main:surgery} $Y_p$ is not surgery on a knot in $S^3$,  
\item \label{main:weight} $\pi_1(Y_p)$ is a weight one group,
\item \label{main:two} $Y_p$ is surgery on a two-component link in $S^3$,
\item \label{main:homcob} no $Y_p$ is smoothly rationally homology cobordant to Auckly's example nor to each other (regardless of orientation). 
\end{enumerate}
\end{theorem}

Theorem~\ref{thm:main} is essentially proved in two steps.  The first step consists of finding an obstruction in Heegaard Floer homology to a homology sphere being surgery on a knot.  The second step consists of an analysis (but not complete computation) of the Heegaard Floer homology of the manifolds $Y_p$.  

Before stating these results, we recall from \cite{OSProperties,OSGraded} that for a homology sphere, its Heegaard Floer homology, $HF^+(Y)$, is a $\Z$-graded $\F[U]$-module, where $\F = \Z/2$ and $U$ lowers degree by 2.  Further, $HF^+(Y)$ admits a non-canonical decomposition  
\[
HF^+(Y) = \Tplus_{d(Y)} \oplus \HFred(Y),
\]
where $\Tplus_{d(Y)}$ is the module $\F[U,U^{-1}]/U \cdot \F[U]$ graded such that $\deg(1) = d(Y)$ and $\HFred(Y)$ is a finite sum of cyclic modules.  The (even) integer $d(Y)$, called the {\em $d$-invariant} or {\em correction term}, is in fact an invariant of smooth rational homology cobordism.  The main obstruction we will present for a homology sphere being surgery on a knot is the following.   

\begin{theorem}\label{thm:obstruction}
Let $Y$ be an oriented integer homology sphere such that $Y = S^3_{1/n}(K)$, for some integer $n$ and some knot $K \subset S^3$.  If $d(Y) \leq -8$, then $U \cdot \HFred_0(Y) \neq 0$.  
\end{theorem}

\begin{remark}
Many others have previously used correction terms to obstruct manifolds from being surgery on a knot (see, for instance, \cite[Corollary 5]{Doig}, \cite[Theorem 4.4]{HW}, \cite[Corollary 9.13, Section 10.2]{OSGraded}).
\end{remark}

\begin{remark}
It is known that $d(S^3) = 0$.  Since Auckly's surgery obstruction required the manifold to be homology cobordant to $S^3$, any manifold one could obstruct from being surgery by Theorem~\ref{thm:obstruction} could not be used for Auckly's argument and vice versa.  

\end{remark}

\begin{remark}
It is straightforward to generalize Theorem~\ref{thm:obstruction} to obtain further restrictions of this form on the Heegaard Floer homology of manifolds with highly negative correction terms obtained by surgery on a knot in $S^3$.  Using such a variant, one can also show that the toroidal Seifert fibered homology sphere $\Sigma(2,5,19,21)$ is not obtained by surgery on a knot.
\end{remark}

In light of Theorem~\ref{thm:obstruction}, we are interested in analyzing both the $d$-invariants of $Y_p$ and the $U$-action on $\HFred(Y_p)$.  

\begin{theorem}\label{thm:computation}
For $p$ a positive, even integer, let $Y_p$ denote the Seifert fibered homology sphere $\Sigma(p,2p-1,2p+1)$, oriented as the boundary of a positive-definite plumbing.  Then, 
\begin{enumerate}[label=(\roman{*}), ref=\roman{*}]
\item \label{computation:d} $d(Y_p) = -p$, 
\item \label{computation:hfred} $U \cdot \HFred_0(Y_p) = 0$. 
\end{enumerate}
\end{theorem}

\begin{figure}[h]
	\includegraphics[width=0.35\textwidth]{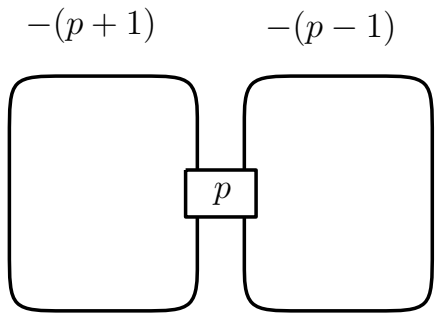}
	\caption{The manifold $Y_p$ presented as surgery on a two component torus link. The box indicates $p$ positive full twists.}
	\label{fig:surgery}
\end{figure}

With this, we may easily prove Theorem~\ref{thm:main}.  
\begin{proof}[Proof of Theorem~\ref{thm:main}]
\eqref{main:surgery}: Notice that the property of a manifold being surgery on a knot in $S^3$ is independent of orientation.  Therefore, we work with $Y_p$ oriented as in Theorem~\ref{thm:computation}.  It is clear that for $p \geq 8$, Theorems~\ref{thm:obstruction} and \ref{thm:computation} now show that $Y_p$ is not surgery on a knot in $S^3$.  

\eqref{main:two}: Since $Y_p$ is a Seifert fibered space with 3 singular fibers, the result follows from \cite[Proposition 8.2]{HW}.  Alternatively, one can directly verify that $Y_p$ is in fact obtained by surgery on the $(2,2p)$ torus link with surgery coefficients $-(p+1)$ and $-(p-1)$ (see Figure \ref{fig:surgery}).   

\eqref{main:homcob}: The result will follow quickly from the following two facts about the rational homology cobordism invariant $d$: $d(-Y) = -d(Y)$ and $d(S^3) = 0$. First we recall that Auckly's examples are homology cobordant to $\Sigma(2,3,5) \# -\Sigma(2,3,5)$ and thus to $S^3$.  Now apply Theorem~\ref{thm:computation}\eqref{computation:d}. 
  
\eqref{main:weight}: This part of the proof was shown to us by Cameron Gordon.  We will show more generally that the Brieskorn sphere $\Sigma(p,q,r)$ has weight one fundamental group.  Suppose that $Z = \Sigma(p,q,r)$ has normalized Seifert invariants $e_0, \frac{p'}{p}, \frac{q'}{q}, \frac{r'}{r}$ (see, for instance, \cite{S}).  Then, we have
\[
\pi_1(Z) = \langle x,y,z,h \mid h \text{ is central}, x^p = h^{p'}, y^q = h^{q'}, z^r = h^{r'}, xyzh^{e_0} = 1 \rangle.  
\]
We claim that $\pi_1(Z)$ is normally generated by $h^{e_0}xy$.  We will show $\pi_1(Z)/ \langle \langle h^{e_0}xy \rangle \rangle$ is trivial.  In this quotient, $z = 1$, so we have 
\[
\pi_1(Z)/ \langle \langle h^{e_0}xy \rangle \rangle \cong \langle x,y,h \mid h \text{ is central}, x^p = h^{p'}, y^q = h^{q'}, h^{r'} = 1, h^{e_0}xy = 1 \rangle. 
\]
Therefore, we can rewrite this as 
\[
\pi_1(Z)/ \langle \langle h^{e_0}xy \rangle \rangle \cong \langle x,h \mid h \text{ is central}, x^p = h^{p'}, (x^{-1}h^{-e_0})^q = h^{q'}, h^{r'} = 1 \rangle.
\]
In particular, $\pi_1(Z) / \langle \langle h^{e_0}xy \rangle \rangle$ is abelian.  However, since $\Sigma(p,q,r)$ is an integer homology sphere, $\pi_1(Z)$ is a perfect group, and thus so is $\pi_1(Z) / \langle \langle h^{e_0}xy \rangle \rangle$.  Therefore, $\pi_1(Z) / \langle \langle h^{e_0}xy \rangle \rangle$ is a perfect abelian group, and thus trivial.  This completes the proof.  
\end{proof}

We are also able to say something for arbitrary homology spheres.  Recall that any reducible homology sphere is not surgery on a knot in $S^3$.  The argument of Gordon and Luecke which is used to prove this result uses that the ambient manifold is $S^3$.  For any homology sphere $Y$, we are able to construct infinitely many reducible manifolds which cannot be surgery on a knot in $Y$.       

\begin{theorem}\label{thm:extra}
Let $Y$ be an integer homology sphere and let $\#_k \Sigma(2,3,5)$ denote the connected sum of $k$ Poincar\'e homology spheres with the same orientation.  For $k \gg 0$, the manifold $\#_k \Sigma(2,3,5)$ is not surgery on a knot in $Y$, regardless of orientation on $Y$.  
\end{theorem} 

\begin{remark}
The reducibility of $\#_k \Sigma(2,3,5)$ is not important for Theorem~\ref{thm:extra}.  What is necessary is a family of integer homology spheres with unbounded $d$-invariants which are L-spaces (i.e., $\HFred = 0$).  The only known irreducible homology sphere L-spaces are $S^3$ and the Poincar\'e homology sphere. 
\end{remark}

{\bf Organization:}  Theorem~\ref{thm:obstruction} is proved in Section~\ref{sec:mappingcone} by utilizing the mapping cone formula for rational surgeries given in \cite{OSRational}.  In Section~\ref{sec:plumbings}, we study the plumbing diagrams of the manifolds $Y_p$ and prove Theorem~\ref{thm:computation}\eqref{computation:d}   using the algorithm of Ozsv\'ath-Szab\'o \cite{OS1}.  In Section~\ref{sec:graded}, we review the algorithm given in \cite{N,CK} to compute the Heegaard Floer homology of Seifert homology spheres.  In Section~\ref{sec:semi}, we analyze $HF^+(Y_p)$ and prove Theorem~\ref{thm:computation}\eqref{computation:hfred}.  Finally, in Section~\ref{sec:extra}, we prove Theorem~\ref{thm:extra}.

\section*{Acknowledgments} We would like to thank Matt Hedden for pointing out that a surgery obstruction could come from comparing the reduced Floer homology with the correction terms.  We would also like to thank Cameron Gordon for supplying the proof of Theorem~\ref{thm:main}\eqref{main:weight}.

%% file: mappingcone.tex
\section{Mapping cones}\label{sec:mappingcone}
The goal of this section is to prove Theorem~\ref{thm:obstruction}.  Let $Y$ be an integer homology sphere with $d(Y) \leq -8$.  Recall that we would like to see that if $Y = S^3_{1/n}(K)$, then $U \cdot \HFred_0(Y) \neq 0$.  We first restrict the possible values of $n$.  
\begin{lemma}\label{lem:negsurgery}
If $Y$ is an integer homology sphere such that $d(Y) < 0$, then $Y$ is not $1/n$-surgery on a knot for any $n < 0$.  
\end{lemma}
\begin{proof}
Suppose that $Y = S^3_{1/n}(K)$, where $n < 0$.  Then, it follows from \cite[Figure 8]{Auckly} that $Y$ is the boundary of a negative-definite four-manifold, $X$.  Since $Y$ is a homology sphere, \cite[Corollary 9.8]{OSGraded} implies $d(Y) \geq 0$.     
\end{proof}

For the rest of this section, we only consider the case of $1/n$-surgery on a knot $K$ for $n>0$.  The main tool is the rational surgery formula of Ozsv\'ath-Szab\'o \cite{OSRational}.  We refer the reader to \cite{NW} for a concise summary.  We very briefly recall the main ingredients for notation without much explanation.  

As usual, let $\Tplus$ denote $\F[U, U^{-1}] / U \cdot \F[U]$.  For each $s \in \mathbb{Z}$, Ozsv\'ath and Szab\'o associate to $K$ a relatively-graded $\F[U]$-module $A_s$, which is isomorphic to the Heegaard Floer homology of a large positive surgery on $K$ in a certain Spin$^c$ structure.  Further, associated to each $s$, there are two graded, module maps $v_s, h_s:A_s \to \Tplus$, which represent maps coming from certain Spin$^c$ cobordisms.  Each $A_s$ admits a splitting $A_s \cong \Tplus \oplus A^{\red}_s$, where $\Tplus$ is the image of $U^N$ for $N \gg 0$ and $A^{\red}_s = \oplus^m_{i=1} \F[U]/U^{k_i}$.  When it will not cause confusion, for $n \geq 0$, we may write $U^{-n}$ to mean the corresponding element of $\Tplus \subset A_s$.  Although $A_s$ is not a module over $\F[U,U^{-1}]$, we will further abuse notation and for an element $a \in \Tplus \subset A_s$, we write $U^{-k}a$ to mean the unique element in $\Tplus \subset A_s$ such that $U^k \cdot U^{-k}a = a$. 

For each $s$, we have that 
\[
v_s|_{\Tplus}(x) = U^{V_s}x
\]
for some non-negative integer $V_s$.  Similarly, 
\[
h_s|_{\Tplus}(x) = U^{H_s}x
\]
for some non-negative integer $H_s$.  Note that each of these maps is surjective.  We will need the following important properties of these integers (see \cite[Section 7]{R}, \cite[Lemma 2.5]{HLZ}, and \cite[Proposition 1.6]{NW}):
\begin{align}
	\label{eqn:symmetry} &H_s =V_{-s}, \\
	\label{eqn:Vineq} &V_s-1 \leq V_{s+1} \leq V_s, \\
	\label{eqn:HV} &H_s = V_s+s, \\
	\label{eqn:gradings} &d(S^3_{1/n}(K)) = -2V_0 = -2H_0.
\end{align}

From this information, we can compute the Heegaard Floer homology of $S^3_{p/q}(K)$ for any rational $p/q \in \mathbb{Q}$.  We will restrict our attention to the case of $S^3_{1/n}(K)$, for $n > 0$.  For each $s$, consider $n$ copies of $A_s$, denoted $A_{s,1},\ldots,A_{s,n}$.  Further, for each $s \in \mathbb{Z}$ and $1 \leq i \leq n$, define $B_{s,i} = \Tplus$.  For an element $x$ in $A_{s,i}$ or $B_{s,i}$, we may write this element as $(x,s,i)$ to keep better track of the indexing.  We will also write $k \pmod{n}$ to refer to the specific representative between 1 and $n$.  Define the map 
\[
\Phi_{1/n} : \bigoplus_{s \in \Z,1 \leq i \leq n} A_{s,i} \to \bigoplus_{s \in \Z, 1 \leq i \leq n} B_{s,i}
\]
by 
\[
\Phi_{1/n}(x,s,i) = (v_s(x),s,i) + (h_s(x),s + \left \lfloor \frac{i}{n} \right \rfloor, i + 1 \hspace{-.08in} \pmod{n}).  
\]
We define an absolute grading on the mapping cone of $\Phi_{1/n}$ (where the $A_{s,i}$ and $B_{s,i}$ are given trivial differential) by requiring that the element $1 \in B_{0,1}$ has grading $-1$ and that $\Phi_{1/n}$ lowers grading by 1.   We remark that the indexing we are using is expressed differently than in \cite{OSRational}.     

\begin{theorem}[Ozsv\'ath-Szab\'o, {\cite[Theorem 1 and Section 7.2]{OSRational}}]\label{thm:rationalsurgery}
The homology of the mapping cone of $\Phi_{1/n}$ is isomorphic to $HF^+(S^3_{1/n}(K))$.  This isomorphism respects the absolute gradings and the $\F[U]$-module structure.      
\end{theorem}

\begin{remark}
Theorem~\ref{thm:rationalsurgery} is not quite stated as in \cite{OSRational}.  Their theorem instead establishes an isomorphism between Heegaard Floer homology and the cone of a chain map whose induced map on homology is $\Phi_{1/n}$.  In general, for a nullhomologous knot in an arbitrary three-manifold, one cannot compute Heegaard Floer homology of surgeries by looking at the cone of the induced map on homology.  However, for knots in $S^3$ (or any L-space), one may compute the homology of the cone of $\Phi_{1/n}$ to obtain the desired result.    
\end{remark}  

With this, we are nearly ready to give the proof of Theorem~\ref{thm:obstruction}.  First, we make an observation about $\HFred$.  Recall that $\HFred(Y)$ is defined to be $HF^+(Y)/ \Im(U^N)$ for $N \gg 0$.  Note that if $a \in HF_{k-2}^+(Y)$ is of the form $U b$ for some $b \in HF_k^+(Y)$ and $a$ is not in $\Im(U^N)$ for $N \gg 0$, then $U \cdot \HFred_k(Y) \neq 0$.  From this and Theorem~\ref{thm:rationalsurgery}, it is straightforward to see that the following proposition implies Theorem~\ref{thm:obstruction}.  

\begin{proposition}
Let $K \subset S^3$ and let $Y=S^3_{1/n}(K)$ for some positive integer $n$. If $d(Y) \leq -8$, then there exist cycles $x$ and $y$ in the cone of $\Phi_{1/n}$ such that
\begin{enumerate}[label=(\roman{*}), ref=\roman{*}]
	\item \label{it:nonzero} $x, y$ are non-zero in homology, 
	\item \label{it:Ux} $y=U x$, 
	\item \label{it:gr0} $\gr(x)=0$, 
	\item \label{it:ynotin} for $N \gg 0$, the element $y$ is not homologous to $U^Nz$ for any cycle $z$.  
\end{enumerate}
\end{proposition}

\begin{proof}
For notation, we let $X$ denote the mapping cone of $\Phi_{1/n}$ and denote by $\mathcal{B}$ the submodule $\bigoplus_{s,1 \leq i \leq n} B_{s, i}$.  Also, let $d=d(Y)$.  Thus $V_0 = H_0 = -\frac{d}{2}$ by \eqref{eqn:gradings}. Since $d \leq -8$, it follows that $V_0 \geq 4$. By \eqref{eqn:Vineq}, we have that $V_1 \geq 3$ and $V_2 \geq 2$.  

We first consider the case when $n=1$.  In this case, we remove the index $i$ used in the $A_{s,i}$ and $B_{s,i}$.  Let $x = U^{1-V_2}$ in $A_2$ and let $y = U x$. Note that $x$ and $y$ are both non-zero in $A_2$ since $V_2 \geq 2$. We have that $v_2(x) = 0$, since $v_2$ restricted to $\Tplus \subset A_2$ is multiplication by $U^{V_2}$. We also have $h_2(x) = 0$, since $h_2$ restricted to $\Tplus \subset A_2$ is multiplication by $U^{H_2}$ and $H_2=V_2+2$ by \eqref{eqn:HV}. Hence, $x$ is a cycle in $X$.  Since $\Phi_{1/n}$ is an $\F[U]$-module map, $y = U x$ must be a cycle in $X$ as well.

We now show that $x$ and $y$ satisfy the conditions of the proposition.    

\begin{enumerate}[label=(\roman{*}), ref=\roman{*}]
	\item Since the image of the differential on $X$ is contained in $\mathcal{B}$ and $x$ is a non-trivial element in $A_2$, the cycle $x$ is non-zero in the homology of $X$. Similarly, $y$ is non-zero in the homology of $X$.
	\item By the definition of $y$, we have that $y = U x$.
	\item Let $z_s$ denote the lowest graded non-zero element of $\Tplus \subset A_s$. Note that $v_s(U^{-V_s}z_s) = h_{s-1}(U^{-H_{s-1}}z_{s-1})$ and this image is the lowest graded non-zero element in $B_s$.  We claim that $\gr(z_0) = -2V_0 = d$.  This follows since $v_0(U^{-V_0}z_0)$ is the lowest graded non-zero element in $B_0$, $v_0$ lowers grading by one, and the grading of the lowest graded non-zero element in $B_0$ is $-1$.  

We have that
	\[ \gr(z_s) = \gr(z_{s-1}) + 2 (H_{s-1} - V_s), \] 
since $v_s(U^{-V_s} z_s) = h_{s-1}(U^{-H_{s-1}} z_{s-1})$, and $v_s$ and $h_{s-1}$ both lower grading by one.

Then
\begin{align*}
	\gr(z_2) &= \gr(z_1) + 2(H_1 - V_2) \\
		&= \gr(z_0) + 2(H_0 - V_1) + 2(H_1 - V_2) \\
		&= \gr(z_0) + 2H_0 +2(H_1-V_1) -2V_2 \\
		&= d - d + 2 -2V_2 \\
		&= 2 -2V_2,
\end{align*}
where the penultimate equality follows from \eqref{eqn:HV} and \eqref{eqn:gradings}.
Since $z_2 = U^{V_2-1} x$ and $U$ lowers grading by two, it follows that $\gr(x) = 0$, as desired.	
	\item We would like to show that for large $N$, $y$ is not homologous to $U^N w$ for any cycle $w$ in $X$. Consider $U^{-N} y \in \Tplus \subset A_2$ for $N \gg 0$. Note that $U^{-N} y = U^{2-V_2-N} z_2 \in A_2$. Then $v_2(U^{-N}y) = U^{2-N} \in B_2 \cong \Tplus$. In particular, $U^{-N} y$ is not in the kernel of the differential on $X$.   Note that 
\[ h_1(U^{2-H_1-N} z_1) = U^{2-N} \in B_2. \]  Moreover, if we choose $N$ greater than
	\[ \max_{s =1, 2 } \left\{ \min \{ n \mid U^n \cdot A^{\red}_s = 0 \} \right\}, \]
then we claim any other element of $X$ whose boundary could cancel with $U^{2-N} \in B_2$ has projection onto $A_1$ given by $U^{2-H_1-N} z_1 \in A_1$.  Indeed, by our choice of $N$, there are no other non-zero elements of $A_1$ or $A_2$ in this grading; further, for an element not contained in $A_1$ or $A_2$, its boundary cannot be contained in $B_2$.  Observe that such an $N$ exists since $A^{\red}_s$ is finite-dimensional as an $\F$-vector space.  
Therefore, if a cycle in $X$ has projection onto $A_2$ given by $U^{-N}y$, then it has projection onto $A_1$ given by $U^{2-H_1-N}z_1$.  

Now, suppose that $w$ is a cycle in $X$ such that $U^N w$ is homologous to $y$.  Then $w$ has projection onto $A_2$ given by $U^{-N}y$.  Thus, the projection of $w$ onto $A_1$ must be $U^{2 - H_1 - N}z_1$. Observe that $2-H_1 < 0$, since $H_1 = V_1 + 1$ and $V_1 \geq 3$ by \eqref{eqn:Vineq}, \eqref{eqn:HV}, and \eqref{eqn:gradings}. Thus, we have that $U^N \cdot U^{2- H_1 - N} z_1 \neq 0$.  This implies that $U^N w$  has non-trivial projection to $A_1$.  Since the image of the differential on $X$ is contained in $\mathcal{B}$ and $y \in A_2$, the cycle $y$ cannot be homologous to an element with non-trivial projection to $A_1$. Hence $y$ is not homologous to $U^N w$.
\end{enumerate}
This completes the proof of the proposition when $n=1$.

The proof when $n>1$ is similar. Let $x = U^{1-V_1} \in \Tplus \subset A_{1, 2}$ and let $y= U x \in A_{1,2}$. As above, it is straightforward to show that $x$ and $y$ are both non-zero in the homology of $X$.  Thus, \eqref{it:nonzero} and \eqref{it:Ux} hold.  We proceed to show that $x$ and $y$ satisfy \eqref{it:gr0} and \eqref{it:ynotin}; the arguments are similar to the $n=1$ case above.
\begin{enumerate} [label=(\roman{*}), ref=\roman{*}]
\setcounter{enumi}{2}
	\item Let $z_{s, i}$ denote the lowest graded element of $\Tplus \subset A_{s, i}$. Note that $\gr(z_{0, 1}) = d$. For $1 \leq i \leq n$,
\begin{equation}\label{eq:zsi} \gr(z_{s, i}) = \gr(z_{s, 1}) + 2s (i - 1),\end{equation}
since 
\begin{itemize}
	\item $h_s(U^{-H_s} z_{s, i}) = v_s(U^{-V_s} z_{s, i+1})$ for $1 \leq i \leq n-1$
	\item $v_s$ and $h_s$ both lower grading by one
	\item $H_s - V_s = s$.
\end{itemize}
We also have that
\begin{equation}\label{eq:z11} \gr(z_{1, 1}) = \gr(z_{0, n}) + 2(H_0 - V_1), \end{equation}
since $v_1(U^{-V_1} z_{1, 1}) = h_0(U^{-H_0} z_{0, n})$. Then, by \eqref{eq:zsi} and \eqref{eq:z11}, 
\begin{align*}
	\gr(z_{1, 2}) &= \gr(z_{1, 1}) + 2 \\
		&= \gr(z_{0, n}) + 2(H_0-V_1) +2 \\
		&= \gr(z_{0, 1}) + 2(H_0-V_1) +2 \\
		&= d - d -2V_1 +2 \\
		&= -2V_1 + 2.
\end{align*}
Since $z_{1, 2} = U^{V_1-1} x$, it follows that $\gr(x) = 0$, as desired.
	\item Consider $U^{-N} y \in \Tplus \subset A_{1, 2}$ for $N \gg 0$. Note that $U^{-N} y = U^{2-V_1-N} z_{1,2} \in A_{1,2}$. Then $v_1(U^{-N} y) = U^{2-N} \in B_{1,2} \cong \Tplus$. In particular, $U^{-N} y$ is not in the kernel of the differential on $X$. Moreover, if we choose $N$ greater than
	\[ \max_{i=1,2 } \left\{ \min \{ n \mid U^n \cdot A^{\red}_{1, i} = 0 \} \right\}, \]
then any other element whose boundary could cancel with $U^{2-N} \in B_{1, 2}$ has projection onto $A_{1,1}$ given by $U^{2-H_1-N} z_{1, 1}$, since
	\[ h_1(U^{2-H_1-N} z_{1,1}) = U^{2-N} \in B_{1,2}. \]
Suppose that $w$ is a cycle in $X$ such that $U^N w$ is homologous to $y$.  Then the projection of $w$ onto $A_{1,1}$ must be $U^{2-H_1-N} z_{1,1}$.  
As discussed above, $H_1 = V_1 + 1$ and $V_1 \geq 3$, and so $2-H_1 <0$.  Thus, $U^N \cdot U^{2 - H_1 - N} z_{1,1} \neq 0$ and therefore, $U^N w$ has non-trivial projection onto $A_{1,1}$. In particular, $y$ is not homologous to $U^N w$, since the image of the differential on $X$ is contained in $\mathcal{B}$.
\end{enumerate}
This completes the proof of the proposition.
\end{proof}

%% file: plumbings.tex
\section{Plumbings}\label{sec:plumbings}
Recall that $Y_p=\Sigma(p,2p-1,2p+1)$, where we have oriented $Y_p$ such that it bounds a positive-definite plumbing.  In this section we determine explicitly the negative-definite plumbing whose boundary is $-Y_p$. We will use this plumbing to compute the correction term of $Y_p$ and hence prove Theorem \ref{thm:computation} \eqref{computation:hfred}.

\begin{figure}[h]
	\includegraphics[width=.40\textwidth]{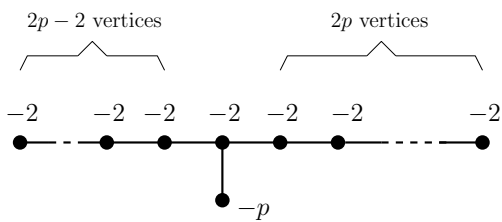}
	\caption{The plumbing graph for $-Y_p$.}
	\label{fig:plumbing}
\end{figure}

\begin{proposition}\label{p:plumb}
For every $p\geq 2$, $-Y_p$ bounds the four-manifold $X_p$ which is the plumbing of disk bundles over spheres intersecting according to the graph in Figure \ref{fig:plumbing}.
\end{proposition}

\begin{proof}
Since $Y_p$ bounds a positive-definite plumbing, clearly $-Y_p$ bounds a negative definite plumbing.  Furthermore, since $Y_p$ has three singular fibers, this plumbing graph will have three arms.  We follow the recipe given in  \cite[Example~1.11]{S} to find this plumbing. We look for the unique integers $e_0,p',q',r'$ solving
$$e_0p(2p-1)(2p+1)+p'(2p-1)(2p+1)+q'p(2p+1)+r'p(2p-1)=-1,$$
with $p>p'\geq 1$, $2p-1>q'\geq 1$, and $2p+1>r'\geq 1$. It can be checked that the solution is given by $e_0=-2$, $p'=1$, $q'=2p-2$, and $r'=2p$. The number $e_0$ is the weight of the central vertex. 

Given integers $m>n>1$, there exists a unique sequence of integers $a_1,\dots,a_k$ with $a_i>1$ for all $i=1,\dots,k$, satisfying
 \begin{equation*}
  \frac{m}{n}= a_1 - \cfrac{1}{a_2
          - \cfrac{1}{\dots
          - \cfrac{1}{a_k} } }.
\end{equation*}
\noindent This is called the \emph{continued fraction expansion} of $\frac{m}{n}$, and we denote it by $\frac{m}{n}=[a_1:a_2:\dots:a_k]$. We now look at the continued fraction expansions of $\frac{p}{p'}$, $\frac{q}{q'}$, and $\frac{r}{r'}$, which determine the negative weights of the vertices on each branch. We have 
\begin{align*}
\frac{p}{1}&=[p], \\
\frac{2p-1}{2p-2}&=[\underbrace{2:2:\dots:2}_{2p-2}],\text{ and}\\
\frac{2p+1}{2p}&=[\underbrace{2:2:\dots:2}_{2p}].
\end{align*}
Hence the result follows.
\end{proof}

We restate Theorem \ref{thm:computation}\eqref{computation:hfred} as the following proposition.  
\begin{proposition} \label{prop:correct}
For $p$ even, we have that $d(Y_p)=-p$.
\end{proposition}
\begin{proof}
Let $X_p$ denote the four-manifold given in Figure \ref{fig:plumbing}. By Proposition \ref{p:plumb} we know $\partial X_p=-Y_p$. A result of Ozsv\'ath and Szab\'o \cite[Corollary 1.5]{OS1} says that the correction term of $-Y_p$ can be computed using the intersection form on $H^2(X_p,\mathbb{Z})$ as follows. Let $\mathrm{Char}(X_p)$ denote the set of all characteristic cohomology classes. Recall that $K\in H^2(X_p,\mathbb{Z})$ is  \emph{characteristic} if $K\cdot [v]+[v]^2\equiv 0\; (\text{mod } 2)$ for every vertex $v$ of the plumbing graph. Next we note that the number of  vertices in the plumbing graph is $4p$. Then the correction term of $-Y_p$ at its unique \Spinc structure is given by
\begin{equation}\label{e:correct}
\displaystyle d(-Y_p)=  \underset{\{K\in\mathrm{Char}(X_p)\}}{\max} \frac{K^2+4p}{4}.
\end{equation}
When $p$ is even, $X_p$ has even intersection form and thus $K=0$ is a characteristic cohomology class. Clearly $K=0$ maximizes the above expression since the intersection form is negative definite. Hence $d(-Y_p)=p$ in this case.  Since $d(-Y_p) = -d(Y_p)$ by \cite[Proposition 4.2]{OSGraded}, we have obtained the desired result.
\end{proof}

\begin{remark} \label{rm:podd}
Though we do not need this for our main argument, we would like to point out that $d(Y_p)=-p+1$ for odd $p$.
\end{remark}

%% file: graded.tex
\section{Graded roots}\label{sec:graded}

The purpose of the present section and the next one is to prove the following result which finishes the proof of Theorem \ref{thm:computation} when combined with Proposition \ref{prop:correct}.
\begin{proposition}\label{prop:rank}
For every even integer $p\geq 4$, we have that $U \cdot \HFred_0(Y_p) = 0$.
\end{proposition}
The proof uses the techniques of graded roots which were introduced by N\'emethi \cite{N} and extensively studied in \cite{CK, KL}.  In this section we motivate and explain our strategy to prove Proposition~\ref{prop:rank} and give the necessary background. The proof will be given in the next section. 

\subsection{Background}\label{subsec:background}
\begin{definition}[N\'emethi, {\cite[Section 3.2]{N}}]
A {\em graded root} is a pair $(\Gamma,\chi)$, where $\Gamma$ is an infinite tree, and $\chi$ is an integer-valued function defined 
on the vertex set  of $\Gamma$ satisfying the following properties. 
\begin{enumerate}[label=(\roman{*}), ref=\roman{*}]
	\item $\chi(u)-\chi(v)=\pm 1$, if there is an edge connecting $u$ and $v$.
	\item $\chi(u)>\mathrm{min}\{\chi (v),\chi (w)\}$, if there are edges connecting $u$ to $v$, and $u$ to $w$.
	\item $\chi$ is bounded below.
	\item $\chi^{-1}(k)$ is a finite set for every $k$.
	\item $|\chi^{-1}(k)|=1$ for $k$ large enough.
\end{enumerate}
\end{definition}

Up to an overall degree shift, every graded root can be described by a finite sequence as follows. Let $\Delta:\{0,\dots,N\}\to\mathbb{Z}$ be a given finite sequence of integers. Let $\tau_\Delta:\{0,\dots,N+1\}\to\mathbb{Z}$ be the unique solution of
\begin{align*}
\tau_\Delta(n+1)-\tau_\Delta(n)=\Delta(n),\text{ with } \tau_\Delta(0)=0.
\end{align*}

\noindent For every $n\in \{0,\dots,N+1\}$, let $R_n$ be the infinite graph with vertex set $\mathbb{Z}\cap [\tau_\Delta (n), \infty)$ and the edge set $\{[k,k+1] \mid k \in \mathbb{Z}\cap [\tau (n), \infty) \}$.  We identify, for each $n\in \{0,\dots,N+1\}$, all common vertices and edges in $R_n$ and $R_{n+1}$ to get an infinite tree $\Gamma_\Delta$. To each vertex $v$ of $\Gamma_\Delta$, we can assign a grading $\chi_\Delta(v)$ which is the unique integer corresponding to $v$ in any $R_n$ to which $v$ belongs.    Clearly many different sequences can give the same graded root. For example the elements $n\in \{0,\dots,N\}$ where $\Delta(n)=0$ do not affect the resulting graded root.

Associated to a graded root $(\Gamma,\chi)$ is a graded $\mathbb{F}[U]$-module $\mathbb{H}(\Gamma)$; we omit the grading function from the notation. As an $\mathbb{F}$-vector space, $\mathbb{H}(\Gamma)$ is generated by the vertices of $\Gamma$.  Further, the grading of a vertex, $v$, is given by $2\chi(v)$.  Finally, $U \cdot v$ is defined to be the sum of vertices $w$ which are connected to $v$ by an edge and satisfy $\chi(w) = \chi(v) - 1$.    

\subsection{Strategy of the proof}
To a large family of plumbed manifolds, N\'emethi associates a graded root whose corresponding module is isomorphic to Heegaard Floer homology up to a grading shift \cite{N}.  In \cite{CK}, N\'emethi's method is simplified for Seifert homology spheres.  Before describing this method in Section~\ref{sec:delta}, we begin with an example to illustrate the process.  This will also enable us to explain the strategy for the proof of Proposition~\ref{prop:rank}.  

For simplicity, we will construct the graded root for $Y_3 = \Sigma(3,5,7)$ and consequently compute its Heegaard Floer homology.  While $Y_3$ does not have $p$ even, this computation will still lend insight into the family of computations we are interested in.  We consider the number $N_{Y_3}=(3\times 5\times 7)-(3\times 5)-(3\times 7)-(5\times 7)=34$. We look at the elements of the semigroup generated by $(3\times 5)$, $(3\times 7)$, $(5\times 7)$, that lie in the interval $[0,N_{Y_3 }]$. The relevant  semigroup elements are $S_{Y_3}=\{0, 15, 21, 30\}$. Let $Q_{Y_3}=\{N_{Y_3}-x \mid x \in S_{Y_3}\}= \{4, 13, 19, 34\}$, and $X_{Y_3}=S_{Y_3}\cup Q_{Y_3}$. Define the function $\Delta_{Y_3}:X_{Y_3}\to \{-1,1\}$ to have value $+1$ on $S_{Y_3}$ and $-1$ on $Q_{Y_3}$. We rewrite the ordered set $X_{Y_3}$, indicating the places where $\Delta_{Y_3}$ is $+1$ with boldface:
$$\{\mathbf{0},4,13,\mathbf{15},19,\mathbf{21}, \mathbf{30},34\}.$$
We write $\Delta_{Y_3}$ as an ordered set, recording in sequence, the value of $\Delta_{Y_3}$ on each element of $X_{Y_3}$:
$$\Delta_{Y_3}=\langle\mathbf{+1},-1,-1,\mathbf{+1},-1,\mathbf{+1},\mathbf{+1},-1\rangle.$$
We then combine the consecutive positive values and the consecutive negative values to write a new sequence which produces the same graded root:
$$\tilde{\Delta}_{Y_3}=\langle \mathbf{+1},-2,\mathbf{+1},-1,\mathbf{+2},-1\rangle.$$

\begin{figure}[h]
	\includegraphics[width=0.40\textwidth]{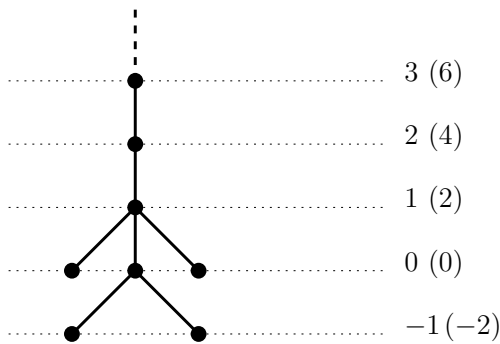}
	\caption{The graded root associated to $Y_3=\Sigma(3,5,7)$.  The grading of each vertex $v$ is written in the form $\chi(v), (\text{absolute grading})$.}
	\label{fig:gradedroot357}
\end{figure}

We indicate the graded root $\Gamma _{Y_3}$ in Figure \ref{fig:gradedroot357}. We can read off the Heegaard Floer homology of $Y_3$ up to a degree shift from its graded root.  As relatively-graded modules, we have that $HF^+(Y_3)$ is isomorphic to 
\[
\mathbb{H}(\Gamma_{Y_3})=\Tplus _{(-2)}\oplus\mathbb{F}_{(-2)}\oplus \mathbb{F}_{(0)}\oplus\mathbb{F}_{(0)}.
\]
Since $d(Y_3) = -2$, we do not shift degrees. Hence, we have 
\[
HF^+(Y_3)=\Tplus _{(-2)}\oplus\mathbb{F}_{(-2)}\oplus \mathbb{F}_{(0)}\oplus\mathbb{F}_{(0)}.
\]

\begin{figure}[h]
	\includegraphics[width=0.50\textwidth]{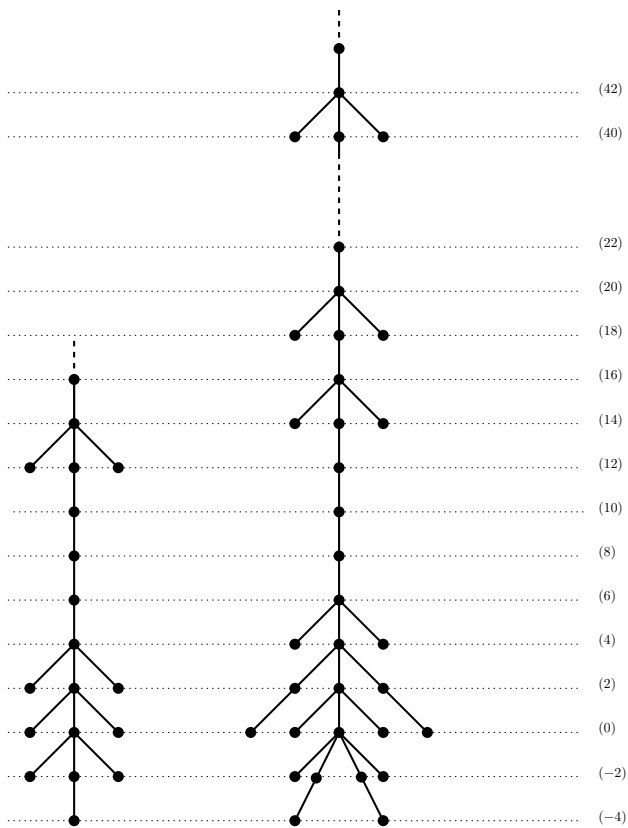}
	\caption{Graded roots associated to $Y_4$ and $Y_5$. Gradings shown correspond to the absolute grading on Heegaard Floer homology.}
	\label{fig:grp45}
\end{figure}

We repeat the same process for $p=4$ and $p=5$. The resulting delta sequences are 
$$\tilde{\Delta}_{Y_4}=\langle 1, -6, 1, -2, 1, -2, 1, -2, 2, -1, 2, -1, 2, -1, 6, -1\rangle,$$
$$\tilde{\Delta}_{Y_5}=\langle 1, -12, 1, -3, 1, -6, 1, -3, 2, -2, 1, -2, 1, -2, 2, -2, 2, -1, 2, -1, 2, -2, 3, -1, 6, -1, 3, -1, 12, -1 \rangle.$$
\noindent See Figure \ref{fig:grp45} for the corresponding graded roots. Since $d(Y_4) = -4$ and $d(Y_5) = -4$, we shift degrees to convert $\mathbb{H}(\Gamma_{Y_p})$ to $HF^+(Y_p)$.

Let us observe why $U \cdot \HFred_0(Y_p)=0$ when $p=3,4,5$ using these graded roots. From the description of the $U$-action on the homology of a graded root $\Gamma$, we see that the dimension of $\ker (U)_n$ is the number of branches ending at degree $n$ whereas $\dim \mathbb{H}_n(\Gamma)$ is the number of vertices in degree $n$. From the pictures of the graded roots of $\Gamma_{Y_p}$ we clearly see that there is exactly one degree $0$ vertex which is not the end of a branch; this vertex is in the image of $U^N$ for $N \gg 0$.  Since $\HFred(Y_p)$ is the cokernel of $U^N$ for $N \gg 0$, we have $U \cdot \HFred_0(Y_p) = 0$ for $p=3,4,5$. 

\begin{figure}[h]
	\includegraphics[width=0.60\textwidth]{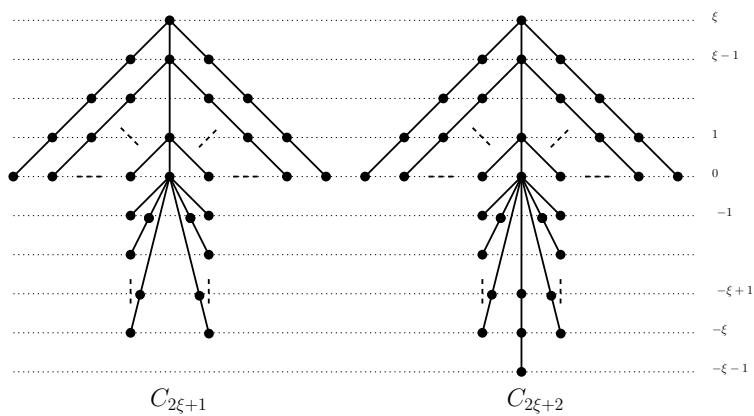}
	\caption{Creatures $\Gamma_{C_p}$ as sub-graded roots of the graded roots associated to $Y_p$. Gradings are $\chi$ values.}
	\label{fig:Creatures}
\end{figure}

In order to prove Proposition \ref{prop:rank} in general, we need to see a pattern in the graded roots of $Y_p$. Repeating the graded root computation for a few more values reveals that the bottom of the graded root of $Y_p$ shows one of the patterns indicated in Figure \ref{fig:Creatures}, depending on the parity of $p$. We call these ``sub-graded roots''  \emph{creatures} and denote them by $\Gamma_{C_p}$.  Proposition \ref{prop:rank} reduces to showing that the bottom of each graded root is the creature $\Gamma_{C_p}$. In order to formalize and prove this we are going to need  abstract delta sequences which were introduced in \cite{KL}.

%% file: delta.tex
\subsection{Delta sequences}\label{sec:delta}  Recall from  \cite{KL} that an \emph{abstract delta sequence} is a pair $(X,\Delta)$, where $X$ is a well-ordered finite set, and $\Delta :X\to \mathbb{Z}-\{0\}$ which is positive at the minimal element of $X$. These objects codify graded roots via the method described in Section~\ref{subsec:background}. 

We review the description of the abstract delta sequence  $(X_Y,\Delta_Y)$ which is associated to an arbitrary  Brieskorn sphere $Y=\Sigma(p,q,r)$. Let $N_Y=pqr-pq-pr-qr$, and let $S_Y$ denote the intersection of the the interval $[0,N_Y]$ with the semigroup generated by $pq, pr$, and  $qr$. Define the set $Q_Y:=\{N_Y-s \mid s \in S_Y \}$, and let $X_Y=S_Y\cup Q_Y$. It turns out that $S_Y$ and $Q_Y$ are disjoint.  Define $\Delta_{Y}:X_Y\to \{-1,1\}$ which takes the value $+1$ on $S_Y$, and $-1$ on $Q_Y$.   Thus, we have $x \in S_Y$ if and only if $N_Y - x \in Q_Y$ and $\Delta_Y(x) = -\Delta_Y(N_Y - x)$ for $x \in X_Y$.  The significance of this abstract delta sequence is the following.
\begin{theorem}[{\cite[Theorem 1.3]{CK}, \cite[Section 11]{N}, and \cite[Theorem 1.2]{OS1}}]\label{thm:nck}
Let $Y = \Sigma(p,q,r)$, oriented as the boundary of a positive-definite plumbing, and let $\Gamma_Y$ be the graded root associated to the abstract delta sequence $(X_Y,\Delta_{Y})$ defined above.  Then, as relatively-graded $\mathbb{F}[U]$-modules, $\mathbb{H}(\Gamma_Y) \cong HF^+(Y)$.
\end{theorem}

\subsection{Refining and merging delta sequences} One can define operations on abstract delta sequences which do not change the corresponding graded root. Two such operations are refinement and merging which we now define. Let $(X,\Delta)$ be a given abstract delta sequence. Suppose there exists a positive integer $t\geq 2$ and an element $z$ of $X$ such that $|\Delta(z)|\geq t$. Pick  integers $n_1,\dots,n_t$, all of which have the same sign as $\Delta(z)$, and satisfy $n_1+\dots+n_t=\Delta (z)$. From this we construct a new delta sequence. Remove $z$ from $X$ and put $t$ consecutive elements $z_1,\dots,z_t$ in its place to get a new ordered set $X'$. Define $\Delta':X'\to \mathbb{Z}$, such that $\Delta '(x)=\Delta(x),\; \text{for all}\;x\in X\setminus \{z\}$, and $\Delta' (z_i)=n_i,\; \text{for all }i=1,\dots,t$. We say that the delta sequence $(X',\Delta')$ is a \emph{refinement} of $(X,\Delta)$ at $z$, and conversely $(X,\Delta)$ is the \emph{merge} of $(X',\Delta')$ at $z_1,\dots,z_t$. 

\begin{definition}\label{d:merref}
An abstract delta sequence is said to be \emph{reduced} if it does not admit any merging (hence there are no consecutive positive or negative values of $\Delta$). An abstract delta sequence is called \emph{expanded}  if it does not admit any refinement (hence every value of delta is $\pm 1$).
\end{definition}

Clearly every delta sequence admits unique reduced and expanded forms. Note that the abstract delta sequences of Brieskorn spheres are in expanded form.

\subsection{Successors and predecessors}\label{subsec:success}
Let $(X,\Delta)$ be an abstract delta sequence.  Denote by $S$ (respectively $Q$) the set of all elements in $X$ where $\Delta$ is positive (respectively negative). For $x \in X$, we define its {\em positive successor} $\suc_+(x) = \min \{x' \in S \mid x < x'\}$ (respectively {\em negative successor} $\suc_-(x) = \min \{x' \in Q \mid x < x'\}$).  In other words, $\suc_{\pm}(x)$ is the first positive/negative element of the delta sequence after $x$.  Should $\suc_{\pm}(x)$ not exist, we treat it as an auxiliary element which is larger than any element in $X$.  Note that $(X,\Delta)$ is in reduced form if and only if for all $x \in S$, $x < \suc_-(x) \leq \suc_+(x)$ and for all $x \in Q$, $x < \suc_+(x) \leq \suc_-(x)$.  We have the analogous notions, $\pre_{\pm}(x)$, which are the {\em predecessors}.

For $x \in S$, define 
\[ \pi_+(x) = \max\{z \in S \mid z < \suc_-(x) \} \quad \textup{ and } \quad \pi_-(x) = \min\{z \in S \mid z > \pre_-(x)\}. \]  
Similarly, for $y \in Q$, we have 
\[ \eta_+(y) = \max \{z \in Q \mid  z < \suc_+(y)\} \quad \textup{ and } \quad \eta_-(y) = \min\{z \in Q \mid z > \pre_+(y)\}. \]   
For $a,b \in X$ with $a<b$, let $[a,b]_X$ denote the set $\{z\in X \mid a\leq z \leq b\}$. Clearly, $[\pi_-(x),\pi_+(x)]_X$ (respectively $[\eta_-(y),\eta_+(y)]_X$) is the maximal set of all consecutive elements of $X$ which are contained in $S$ (respectively $Q$) which contains $x$ (respectively $y$).  

We describe an explicit model for the reduced form of $(X,\Delta)$, denoted $(\tilde{X}, \tilde{\Delta})$ such that $\tilde{X} \subset X$.   This is done as follows.  Define $\tilde{S} = \{\pi_+(x) \mid x \in S\}$ (i.e., the largest endpoints of each maximal interval of elements with positive values) and $\tilde{Q} = \{\eta_-(y) \mid y \in S\}$.  We then merge each $x \in [\pi_-(x),\pi_+(x)]_X$ (respectively $y \in [\eta_-(y),\eta_+(y)]_X$) with $\pi_+(x)$ (respectively $\eta_-(y)$).  Consequently, if $\Delta$ is expanded, we have that $\tilde{\Delta}(\pi_+(x)) = \# [\pi_-(x),\pi_+(x)]_X$ and $\tilde{\Delta}(\eta_-(y)) = - \# [\eta_-(y),\eta_+(y)]_X$.  

When discussing the reduced form of $(X,\Delta)$, we will always assume we are working with this explicit model for the reduced form of $(X,\Delta)$. There is also an obvious quotient from $X$ to $\tilde{X}$ given by $x \mapsto \pi_+(x)$ for $x \in S$ and $y \mapsto \eta_-(y)$ for $y \in Q$.  We will sometimes not distinguish between an element  of  $X$ and its representative in $\tilde{X}$.  The reason for this is that if $x < y$ or $y < x$ in $X$, for $x \in S$ and $y \in Q$, then the same inequality holds for their images in $\tilde{X}$.   

Let  $Y$ be a Brieskorn sphere. We would like to study the reduced form $(\tilde{X}_Y,\tilde{\Delta}_Y)$ of the delta sequence $(X_Y,\Delta_Y)$ according to the above model.  First we point out that $N_Y-x\in \tilde{X}_Y$ whenever $x\in \tilde{X}_Y$. Moreover, it follows from the construction that $\tilde{\Delta}_Y(N_Y-x) = -\tilde{\Delta}_Y(x)$ for $x \in \tilde{X}_Y$. Finally we have the following. 

\begin{proposition} For a Brieskorn sphere $Y$, we have that $0$ and $N_Y$ are both contained in $\tilde{X}_Y$.  
\end{proposition}

\begin{proof}
Let $Y=\Sigma(p,q,r)$ be a Brieskorn sphere. Let $x_0$ denote the minimum of $\{pq,pr,qr\}$.  Throughout this proof, all successors and predecessors are taken with respect to $X_Y$ and not $\tilde{X}_Y$.  Note that $\suc_+(0)=x_0$.  Let $y$ denote the maximal element of  $S(pq,pr,qr)$ less than $N_Y$.  Then $\suc_-(0)=N_Y-y$. Since $x_0$ is the smallest generator of $S(pq,pr,qr)$, gaps between any two elements of this semigroup are less than or equal to $x_0$. Since $N_Y\not \in S(pq,pr,qr)$, we have $N_Y-y<x_0$ by definition of $y$. Hence, we have $0<\suc_-(0)<\suc_+(0)$ which implies $\pi_+(0)=0$. Therefore $0\in \tilde{X}_Y$. 
Consequently, $N_Y-0 = N_Y \in \tilde{X}_Y$.  
\end{proof}

\subsection{Tau functions and sinking delta sequences} Given an abstract delta sequence $(X,\Delta)$, one defines the well-ordered set $X^+:=X\cup \{z^+ \}$ where $z^+>z$ for all $z\in X$, and a function $\tau_\Delta:X^+\to \mathbb{Z}$, as in Section~\ref{subsec:background}, with the following formula
\begin{equation*}
\tau_\Delta(z)=\underset{w<z}{\sum_{w\in X}}\Delta(w),\;\mathrm{for}\;\mathrm{all}\;z\in X^+.
\end{equation*} 
\noindent We call $\tau_\Delta$ the \emph{tau function} of the delta sequence $(X,\Delta)$.  An important part of the study of abstract delta sequences is to detect where their tau functions attain their absolute minimum.  Below we define a class of delta sequences whose tau functions have easily detectable minimum.
\begin{definition}\label{d:sink}
Let $(X,\Delta)$ be an abstract delta sequence and let $(\tilde X, \tilde \Delta)$  be its reduced form. We say that $(X,\Delta)$ is \emph{sinking} if 

\begin{enumerate}[label=(\roman{*}), ref=\roman{*}]
\item \label{i:sink1} the maximal  element $z_{\max}$ of $X$ belongs to $Q$, 
\item \label{i:sink2} for  every element $x\in \tilde{S}$, we have $\tilde \Delta(x)\leq |\tilde \Delta(\suc_-(x))|$,  and
\item \label{i:sink3} $\tilde \Delta(\pre_+(z_{\max}))< |\tilde \Delta (z_{\max})|$.
\end{enumerate}

\end{definition}
\begin{proposition}\label{p:sink}
The tau function of a sinking delta sequence attains its absolute minimum at the last element and nowhere else.
\end{proposition}
\begin{proof}
Follows from the definition.
\end{proof}
We will also need certain dimension formulas for $\mathbb{H}(\Gamma_\Delta)$, regardless of whether $(X,\Delta)$ is sinking. We find it convenient to work in the reduced form. Let $(X,\Delta)$ be a given abstract delta sequence and let $(\tilde{X},\tilde{\Delta})$ be its reduced form. Let $\tilde{\tau}: \tilde{X}^+\to \mathbb{Z}$ be the tau function of the reduced sequence. For any $z \in \tilde{X}^+$ other than the minimal element, let $\pre(z)$ denote the immediate predecessor of $z$ in $\tilde{X}^+$ (i.e., $\pre_-(z)$ if $z \in \tilde{S}$ and $\pre_+(z)$ if $z \in \tilde{Q}$). Denote by $z_{\min}$ the minimal element of $\tilde{X}^+$.
\begin{proposition}\label{p:rankform}
 We have
\begin{enumerate}[label=(\roman{*}), ref=\roman{*}]
\item \label{i1:rankform} $\dim (\mathbb{H}_{0}(\Gamma_\Delta))= \# \{z \in \tilde{X}^+\setminus \{z_{\min}\} \mid \tilde{\tau} (\pre(z))>0 \text{ and } \tilde {\tau} (z)\leq 0 \} +1$,
\item \label{i2:rankform} $\dim (\ker U)_{0}=\# \{z \in \tilde{X}^+ \setminus \{z_{\min}\} \mid \tilde{\tau} (\pre(z))>0 \text{ and }  \tilde{\tau} (z)=0\} +1$.
\end{enumerate}
\end{proposition}
\begin{proof}
As we noted before, passing to the reduced form does not change the graded root. The proof follows from the construction of the graded root $\Gamma_\Delta$ and the description of the $U$-action on $\mathbb{H}(\Gamma_\Delta)$.

First of all $\dim (\mathbb{H}_{0}(\Gamma_\Delta))$ is the total number of vertices with degree $0$. Now recall the construction of the graded root from the tau function given in Section~\ref{subsec:background}. Since $\tilde \Delta$  is reduced, each time   $\tilde{\tau} (\pre(z))>0 \text{ and } \tilde {\tau} (z)\leq 0$ happens, we create a new vertex of degree $0$. Since $\tilde{\tau}(z_{\min})=0$, we have one more vertex of degree $0$.   

By the description of this $U$-action, $\dim (\ker U)_{0}$ equals the number of valency one vertices on $\Gamma_\Delta$ which have degree $0$.  Each time $\tilde{\tau} (\pre(z))>0 \text{ and }  \tilde{\tau} (z)=0$, happens we create such a vertex. We need to add one to include the vertex corresponding to $z_{\min}$. 
\end{proof}
\subsection{Symmetric delta sequences} We see an obvious symmetry in Figures \ref{fig:gradedroot357} and \ref{fig:grp45}. In fact this symmetry more generally holds for the graded roots of Seifert homology spheres. The purpose of the next two definitions is to characterize those delta sequences whose graded roots are symmetric. To simplify the description, we shall use the notation $f=\langle k_1,k_2,\dots,k_n\rangle$ to denote the function $f:X\to \mathbb{Z}$, whose value is $k_1$ at the minimal element of $X$, then $k_2$ at the successor of the minimal element, et cetera.
\begin{definition}\label{d:ope}
Given an abstract delta sequence $\Delta=\langle k_1,\dots,k_n\rangle$, we define the following
\begin{enumerate}[label=(\roman{*}), ref=\roman{*}]
\item \emph{negation}, $-\Delta =\langle -k_1,\dots,-k_n\rangle$,
\item \emph{reverse}, $\overline{\Delta}=\langle k_n,\dots,k_1\rangle$. 
\end{enumerate}
Note that neither the negation nor the reverse of an abstract delta sequence need be an abstract delta sequence.  If $\Delta_1=\langle k_1,\dots,k_n\rangle $ and $\Delta_2=\langle \ell_1,\dots,\ell_m\rangle $ are abstract delta sequences we define their \emph{join} by 
$$\Delta_1 * \Delta_2= \langle k_1,\dots,k_n, \ell_1,\dots,\ell_m\rangle.  $$
\end{definition}
\begin{definition}\label{d:sym}
Let 
$$\Delta=\langle k_1,\dots,k_n\rangle $$
be an abstract delta sequence. The {\em symmetrization} of $\Delta$ is the following abstract delta sequence
$$\Delta^{\mathrm{Sym}}:=\Delta * - \overline{\Delta}=\langle k_1,\dots,k_n,-k_n,\dots, -k_1\rangle .$$
\end{definition}
\begin{remark}\label{r:sym}
Since $x\in S_Y$ if and only $N_Y-x\in Q_Y$ and $\Delta_Y(x)=-\Delta_Y(N_Y-x)$, we see that $N_Y/2 \not \in X_Y$; hence we have
$$\Delta_Y= (\Delta_Y|_{[0,N_Y/2]})^\textrm{Sym}.$$
Further, observe that by the symmetry of $\Delta_Y$, if the maximal element of $X_Y \cap [0,N_Y/2]$ is an element of $S_Y$ (respectively $Q_Y$), then the minimal element of $X_Y \cap [N_Y/2,N_Y]$ is an element of $Q_Y$ (respectively $S_Y$).  Therefore, it follows that $\tilde{\Delta}_Y = (\tilde{\Delta}_Y|_{[0,N_Y/2]})^\textrm{Sym}$.
\end{remark}

%% file: semigroups.tex
\section{Semigroups and creatures}\label{sec:semi} 
Having given the necessary background, we are now ready to prove Proposition \ref{prop:rank}. First we formally define the creatures given in Figure \ref{fig:Creatures}, by indicating their delta sequences.  Then we will observe that Proposition \ref{prop:rank} holds for the creature graded roots, which will be denoted $\Gamma_{C_p}$; namely, we will show $U \cdot \mathbb{H}^{\textup{red}}_0(\Gamma_{C_p}) = 0$. Finally we shall prove a technical decomposition theorem which essentially reduces the proof of Proposition \ref{prop:rank} for $Y_p$ to checking that it holds for the creatures.  Throughout this section, let $p$ be an even integer with $p \geq 4$.  We will often write $p = 2\xi+2$.  
\begin{definition}\label{d:cre}
For every $p=2\xi+2$ with $\xi\geq 1$, the \emph{creature} $\Gamma_{C_p}$ is the graded root defined by the symmetrization of the abstract delta sequence
$$\Delta_{C_p}=\langle \xi,-\xi,(\xi-1),-(\xi-1),\dots, 2,-2,1,-2,1,-2,2,\dots,-(\xi-1),(\xi-1),-\xi,\xi,-(\xi+1)\rangle .$$
\end{definition}

Let $p$ be given, and consider the creature graded root $\Gamma_{C_p}$ and its homology $\mathbb{H}(\Gamma_{C_p})$, which is an $\mathbb{F}[U]$-module 
supported in even degrees (see Section~\ref{subsec:background} for the construction of $\mathbb{H}(\Gamma)$).  

\begin{proposition}\label{p:crea} For $p=2\xi+2$ with $\xi\geq 1$, we have that
$$\dim(\ker U)_0 +1 = \dim \mathbb{H}_{0}(\Gamma_{C_p}).$$
\end{proposition}
\begin{proof}
Note that $(\Delta_{C_p})^{\mathrm{Sym}}$ is in reduced form and its tau function is given by 
\begin{align} \label{e:taucre}
\tau = \langle 0,\xi,0,\xi-1,0,\xi-2,0,\dots,2,0,1,\underline{-1},0,-2,0,\dots,-(\xi-1),0,-\xi,0,-(\xi+1),\\
\nonumber  0,-\xi,0,-(\xi-1),0,\dots,-2,0,-1,1,0,2,0,\dots,\xi-1,0,\xi,0 \rangle.
\end{align}
Let $\pre(z)$ denote the immediate predecessor of $z$, as in Proposition~\ref{p:rankform}. From this we observe that there is exactly one element that belongs to the set $\{z \mid \tau(\pre(z))>0 \text{ and } \tau(z) \leq 0\}$,   but does not belong to the set $\{z \mid \tau(\pre(z))>0 \text{ and } \tau(z)=0\}$. We underlined the value of the tau function at this element in \eqref{e:taucre} above. Hence by Proposition \ref{p:rankform}, the proof follows.
\end{proof}

Let $Y_p=\Sigma(p,2p-1,2p+1)$. Let $\Delta_{Y_p}$ denote the corresponding abstract delta sequence as described in Section~\ref{sec:delta}, and let $\tilde \Delta_{Y_p}$ denote its reduced form. The proof of Proposition  \ref{prop:rank} will follow from the following technical statement about $\tilde \Delta_{Y_p}$.

\begin{lemma}\label{l:dec}
For every even integer $p\geq 4$, we have the following decomposition
$$\tilde \Delta_{Y_p}=(\Delta_{Z_p}*\Delta_{C_p})^{\mathrm{Sym}},$$
where $\Delta_{C_p}$ is the creature sequence given in Definition \ref{d:cre}, and  $\Delta_{Z_p}$ is a sinking delta sequence.
\end{lemma}

\noindent Let us first see why the above lemma implies Proposition~\ref{prop:rank}, and thus completes the proof of Theorem~\ref{thm:computation}.

\begin{proof}[Proof of Proposition \ref{prop:rank}]

Consider the graded root $\Gamma_{Y_p}$ whose grading is shifted so that it agrees with the absolute grading of $HF^+(Y_p)$ (see Theorem~\ref{thm:nck}). The decomposition in Lemma~\ref{l:dec} implies that the creature graded root $\Gamma_{C_p}$ embeds into $\Gamma_{Y_p}$ as a subgraph. Moreover, Proposition~\ref{prop:correct}  implies that this embedding is in fact degree preserving. Since $\Delta_{Z_p}$ is sinking, by Proposition~\ref{p:sink} and \eqref{e:taucre}, the minimum value of $\tau_{\Delta_{Z_p}}$ is 0, which is the initial value of $\tau_{\Delta_{C_p}}$, and is uniquely attained at the maximal element of $\Delta_{Z_p}$.  Thus, we see that $\mathbb{H}_{\leq 0}(\Gamma_{Y_p}) = \mathbb{H}_{\leq 0}(\Gamma_{C_p})$ as graded $\F[U]$-modules.  By Theorem~\ref{thm:nck} and Proposition \ref{p:crea}, we see $\dim HF^+_0(Y_p) = \dim(\ker U)_0 + 1$.  This implies the desired result.
\end{proof}

The proof of Lemma \ref{l:dec} occupies the rest of the section. Let $r_\pm = p(2p \pm 1)$ and $w = (2p-1)(2p+1)$.  Therefore, in order to study the abstract delta sequence for $\Sigma(p,2p-1,2p+1)$, we must work with the semigroup $S(r_-,r_+,w)$ generated by $r_-$, $r_+$ and $w$.   Ideally one would like to  describe explicitly the elements of $S(r_-,r_+,w)\cap[0,N_{Y_p}]$. This set seems to be too complicated  at the moment. Instead we will only need an explicit description of $S(r_-,r_+,w)\cap [0,(p-1)r_+]$; note that $(p-1)r_+<N_{Y_p}$.  We begin with an important subset of $S(r_-,r_+,w) \cap [0,(p-1)r_+]$.   

\begin{lemma}\label{lem:sem2}
The intersection $S(r_-,r_+) \cap [0,(p-1)r_+]$, as an ordered set, is given by the ordered set
\begin{align*}
\{&0,\\
 &r_-,\; r_+,\\
 &2r_-,\; r_- + r_+,\; 2r_+,\\
 &3r_-,\; 2r_- + r_+,\; r_- + 2r_+,\; 3r_+,\\
 &\; \vdots\\
 &(p-1)r_-,\; (p-2)r_- + r_+,\ldots,\; (p-1)r_+\}.  
\end{align*}
\noindent Here we deliberately break the lines, so the pattern of the  elements  is visible.
\end{lemma}
\begin{proof}
First, it is clear that $(p-1)r_+ \in S(r_-,r_+)$, and therefore is the maximal element.  Next, note that if $a + b = k$ and $1 \leq a \leq k$ and $1 \leq b \leq k$, then since $r_-<r_+$, 
\[
a r_- + (b-1) r_+ < (a-1) r_- + b r_+.  
\]
Therefore, to establish the order as given in the statement of the lemma, we just need to show that as long as $k \leq (p-1)$, we have $k r_+ < (k+1) r_-$.  This inequality is easily checked using the definition $r_{\pm} = p(2p \pm 1)$.   
\end{proof}

\begin{fact}\label{fact:consec}
Let $a$ and $b$ be positive integers. Observe that since $w = r_- + r_+ - 1$, the sequence of elements of $S(r_-,r_+,w)$  
\begin{align*}
&(a - \min\{a,b\})r_- + (b - \min\{a,b\}) r_+ + \min\{a,b\} w, \\
&(a - \min\{a,b\} + 1)r_- + (b-\min\{a,b\} + 1) r_+ + (\min\{a,b\} -1) w, \\ 
&\; \vdots\\
&(a-1)r_- + (b-1)r_+ + w, \\
&ar_- + b r_+
\end{align*}
is consecutive in $\mathbb{N}$.  This will be used frequently throughout the proof.  
\end{fact}

Before proceeding, we point out that if $x \in S_{Y_p}$ can be written as $x= ar_- + br_+ + cw$ for some non-negative integers $a,b$, and $c$, then this decomposition is unique by the Chinese remainder theorem.   Suppose now that $a$ is a non-negative integer and $b$ is a  positive integer such that $a+b \leq p-1$. Fact~\ref{fact:consec}, combined with this observation about the unique representability of elements in $S_{Y_p}$, implies  
\begin{align}\label{eq:consec}
\nonumber (a+1) r_- + (b-1) r_+ < &(a - \min\{a,b\})r_- + (b - \min\{a,b\}) r_+ + \min\{a,b\} w\\
& \vdots  \\
\nonumber < &(a-1)r_- + (b-1)r_+ + w,\\
\nonumber < & ar_-+ br_+.  
\end{align}

Combining Lemma~\ref{lem:sem2} with \eqref{eq:consec} now gives a complete description of $S(r_-,r_+,w) \cap [0,(p-1)r_+]$.  Namely for each element $x=ar_-+br_+$ of the pyramid given in Lemma~\ref{lem:sem2}, we have $\min\{a,b\}$ more consecutive elements preceding $x$. 
In particular, there are no elements of $S(r_-,r_+,w)$ between $(k-1)r_+$ and $kr_-$ for $k \leq p-1$.  

We also point out two inequalities which will be used throughout the proof of Lemma~\ref{l:dec}:  
\begin{align}
\label{i:fre1} (p-1)r_-+(p-3)r_+ &<N_{Y_p},\\
\label{i:fre2} (p-2)r_-+(p-2)r_+&>N_{Y_p}.
\end{align}
The validity of these two inequalities can be checked directly from the definitions $r_\pm=p(2p\pm 1)$ and $N_{Y_p}=4p^3-8p^2-p+1$.

\subsection{The proof of Lemma~\ref{l:dec}}

In order to study $\Delta_{Y_p}$ we will find it more convenient to work with its reduced form.   We shall make use of the explicit model of the reduced form given in Section~\ref{subsec:success}. Hence, we will heavily rely on the notation introduced there.  Recall that in order to determine the reduced form, we must compute $\pi_{\pm}(x)$ for $x \in S_{Y_p}$.  

\begin{lemma}\label{lem:intervals}
Suppose that $x \in S_{Y_p}$ satisfies $x = a r_- + br_+$, where $a,b \geq 0$, and that $x \leq 2r_- + (p-3)r_+$.  Then, 
\begin{enumerate}[label=(\roman{*}), ref=\roman{*}]
 \item  \label{shufflei:1} $x < N_{Y_p} - (p-a-1) r_- - (p-b-3)r_+ < \suc_+(x)$ 
 \item \label{shufflei:2} $[\pi_-(x),\pi_+(x)]\cap S_{Y_p} = \{x-\min\{a,b\},\ldots,x\}$, unless $x = (p-2)r_+$ or $(p-1)r_-$.  In either of these exceptional cases, we have $[\pi_-(x),\pi_+(x)] \cap S_{Y_p} = \{(p-2)r_+, (p-1)r_-\}$.   
\end{enumerate}  
\end{lemma} 
\begin{proof}
First, let 
$$x' = \left \{
\begin{array}{ll}
(a-1) r_- + (b+1)r_+ & \text{ if } a \geq 1\\
(b+1)r_-&  \text{ if } a = 0.
\end{array} \right .$$  
One can check that $x' \leq (p-1)r_+$. Hence we have by Lemma~\ref{lem:sem2} and Fact~\ref{fact:consec} that $x < x'$ and that the elements of $S_{Y_p}$ strictly between $x$ and $x'$ are of the form $x' - i$, for $1 \leq i \leq \min\{a-1,b+1\}$ and they are consecutive in $X_{Y_p}$.  Thus, 
$$\suc_+(x) = \left \{
\begin{array}{ll}
x' - \min\{a-1,b+1\} & \text{ if } a \geq 1\\
x'&  \text{ if } a = 0.
\end{array} \right .$$  
Since the elements between $\suc_+(x)$ and $x'$ are consecutive in $X_{Y_p}$ (and $S_{Y_p}$), if $y \in Q_{Y_p}$ satisfies $y < x'$, it must satisfy $y < \suc_+(x)$.  However, it is straightforward to verify that \eqref{i:fre1} and \eqref{i:fre2} imply 
\begin{equation}\label{e:inbetween}
x < N_{Y_p} - (p-a-1) r_- - (p-b-3)r_+ < x'. 
\end{equation}
This completes the first part of the lemma.  

In fact, when $b < (p-2)$, then $p - a - 1 \geq 0$ and $p - b - 3 \geq 0$, so we have found an element of $Q_{Y_p}$ between $x$ and $\mathrm{suc}_+(x)$, and thus $\pi_+(x) = x$.  When $b \geq (p-2)$, since $x \leq 2r_- + (p-3)r_+$, we have $x = (p-2)r_+$ by Lemma~\ref{lem:sem2}.  By \eqref{e:inbetween}, we do see $(p-2)r_+ < N_{Y_p} - (p-1)r_- + r_+ < x' = (p-1)r_-$, but it turns out that $N_{Y_p} - (p-1)r_- + r_+$ is not of the form $N_{Y_p} - z$, for any $z \in S_{Y_p}$.  We will deal with this exceptional case shortly.    

First, we would like to determine $\pi_-(x)$ in the generic case. By Fact \ref{fact:consec}, $\{x - \min\{a,b\},\ldots,x\}$ is a consecutive subset of $X_{Y_p}$ which is contained in $S_{Y_p}$. Let
$$x'' = \left \{
\begin{array}{ll}
(a+1) r_- + (b-1)r_+ & \text{ if } b \geq 1\\
(a-1)r_+&  \text{ if } b = 0.
\end{array} \right .$$  
\noindent By \eqref{eq:consec} and Lemma \ref{lem:sem2}, we have that $\pre_+(x - \min\{a,b\})=x''$. Again, Lemma \ref{lem:sem2}, \eqref{i:fre1}, and \eqref{i:fre2} imply 
\[
x'' < N_{Y_p} - (p-a-2) r_- - (p-b-2)r_+ < x. 
\]
\noindent Similar to above, when $a < p-1$, we have that $N_{Y_p} - (p-a-2)r_- - (p-b-2)r_+ \in Q_{Y_p}$, since $p - a -2, p -b-2 \geq 0$.  In this case, $\pi_-(x) = x - \min\{a,b\}$.  Thus, if $x$ is neither $(p-2)r_+$ nor $(p-1)r_-$, the second claim follows, since $\pi_+(x) = x$.  

In order to deal with the exceptional cases, we will prove that there is no element of $Q_{Y_p}$ between $(p-2)r_+$ and $(p-1)r_-$. Since the above arguments show that there exists an element of $Q_{Y_p}$ between $\pre_+((p-2)r_+)$ and $(p-2)r_+$, and an element of $Q_{Y_p}$ between $(p-1)r_-$ and $\suc_+((p-1)r_-)$, this will establish that $[\pi_-(x),\pi_+(x)] \cap S_{Y_p} = \{(p-2)r_+,(p-1)r_-\}$ for $x = (p-2)r_+$ or $(p-1)r_-$.  Here, we are using our description of $S_{Y_p}$ to deduce that there are no elements of $S_{Y_p}$ between $(p-2)r_+$ and $(p-1)r_-$.  

Suppose $y \in Q_{Y_p}$ satisfies $(p-2)r_+ < y < (p-1)r_-$.  Then, write $y = N_{Y_p} - z$, where $z \in S_{Y_p}$.  We therefore have $(p-2)r_+ + z < N_{Y_p}$ and $(p-1)r_- + z > N_{Y_p}$.  By \eqref{i:fre1} and \eqref{i:fre2}, we have $z < (p-2)r_-$ and $z > (p-3)r_+$.  However, there is no element of $S_{Y_p}$ between $(p-3)r_+$ and $(p-2)r_-$.  This is a contradiction.  Thus, there are no elements in $Q_{Y_p}$ between $(p-2)r_+$ and $(p-1)r_-$, which is what we needed to show.  
\end{proof}

We remark that more generally, if $x = ar_- + br_+$ and $x \geq (p-1)r_+$, we are still able to deduce that $\{x - \min\{a,b\},\ldots,x \} \subset [\pi _-(x),\pi_+(x)]\cap S_{Y_p}$.  Finally, recall that for $x \in X_{Y_p}$, we may also write $x$ for the induced element in $\tilde{X}_{Y_p}$.  

\begin{proposition}\label{prop:reducedstructure}
The reduced form $\tilde{\Delta}_{Y_p}$ of $\Delta_{Y_p}$ has the following properties:
\begin{enumerate}[label=(\roman{*}), ref=\roman{*}]
\item \label{i1:reducedstructure} As ordered subsets of $\mathbb{N}$, $\tilde{S}_{Y_p} \cap [0,2r_- + (p-3)r_+] = S(r_-,r_+)  \cap [0,2r_- + (p-3)r_+]\setminus\{(p-2)r_+ \}$. 
\item \label{i2:reducedstructure} Let $x \in S(r_-,r_+) \cap [0,2r_- + (p-3)r_+]\setminus\{(p-2)r_+, (p-1)r_- \}$ be expressed as $x = ar_- + br_+$.  Then $\tilde{\Delta}_{Y_p}(x) = \min\{a,b\} + 1$. We also have $\tilde{\Delta}_{Y_p}((p-1)r_-)=2$.
\item \label{i3:reducedstructure} Let $x \in \tilde{S}_{Y_p}$ and suppose $x < N_{Y_p} - cr_- - dr_+ < \suc_+(x)$, where $c,d \geq 0$.  Then $\tilde{\Delta}_{Y_p} (\suc_-(x)) \leq -\min\{c,d\}-1$.  
\end{enumerate}    
\end{proposition}
\begin{proof}
Recall the construction of $\tilde{S}_{Y_p}$ given in Section~\ref{subsec:success}, namely that $\tilde{S}_{Y_p}$ consists of the elements of the form $\pi_+(x)$ for $x \in S_{Y_p}$.  For this proof, all predecessors and successors will be taken with respect to $X_{Y_p}$ and never $\tilde{X}_{Y_p}$.  Let $x \in S_{Y_p} \cap [0,2r_- + (p-3)r_+]\setminus\{(p-2)r_+, (p-1)r_- \}$ and suppose that $x$ is of the form $x = ar_- + br_+$.  Lemma~\ref{lem:intervals}\eqref{shufflei:2} gives that if $x \neq (p-2)r_+, (p-1)r_-$, then $[\pi_-(x),\pi_+(x)]\cap S_{Y_p} = \{x-\min\{a,b\}, \ldots,x\}$ and thus $\pi_+(x) = x$. It is also shown in Lemma~\ref{lem:intervals}\eqref{shufflei:2} that $\pi_+((p-2)r_+) = \pi_+((p-1)r_-) =(p-1)r_-$, so $(p-2)r_+$ must be excluded from $\tilde{S}_{Y_p}$. This finishes the proof of the first claim.

Since $\Delta_{Y_p}$ is expanded, 
\[
\tilde{\Delta}_{Y_p}(x) = |\{x-\min\{a,b\},\ldots,x\}| = \min\{a,b\}+1,
\] 

\noindent for every $x\in \tilde{S}_{Y_p} \cap [0,r_- + (p-2)r_+]\setminus \{ (p-1)r_-\}$. In the exceptional case where  $x=(p-1)r_-$, we have $\pi_-(x)=(p-2)r_+$ and $\pi_+(x)=(p-1)r_-$, hence $\tilde{\Delta}_{Y_p}(x)=|\{(p-2)r_+, (p-1)r_-\}|=2$.  This proves the second claim.

It remains to establish the final claim in the proposition. Let $y  = N_{Y_p} - cr_- - dr_+$ and suppose that $x < y < \suc_+(x)$.  Then, we must have $\suc_-(x) \in [\eta_-(y),\eta_+(y)]$.  This implies that $\tilde{\Delta}_{Y_p}(\suc_-(x)) = - \# [\eta_-(y), \eta_+(y)]\cap X_{Y_p}$.  As discussed above, for any $z = s r_- + t r_+ \in S_{Y_p}$, regardless of whether $s + t \leq p-1$, we have that $\{z-\min\{s,t\},\ldots,z\} \subset [\pi_-(z),\pi_+(z)]$.  From this, we can deduce that  
\[
\{N_{Y_p} - cr_- - dr_+, \ldots, N_{Y_p} - cr_- - dr_+ + \min\{c,d\} \} \subset [\eta_-(y),\eta_+(y)].
\]
Therefore, we must have $\tilde{\Delta}_{Y_p}(\suc_-(x)) \leq -\min\{c,d\} - 1$.  This completes the proof.  
\end{proof}

In order to prove Lemma~\ref{l:dec} we are interested in finding a decomposition $\tilde \Delta_{Y_p} = (\Delta_{Z_p} * \Delta_{C_p})^{\Sym}$, such that $\Delta_{Z_p}$ is sinking and $\Delta_{C_p}$ is the creature sequence from Definition~\ref{d:cre}.  
Recall that we write $p=2\xi+2$ for some positive integer $\xi$. Define 
\begin{equation}\label{e:defK}
K  = (\xi-1) r_- + \xi r_+.
\end{equation} 
\noindent Note that 
\begin{equation}\label{e:Kest}
K < (p-3) r_+ < N_{Y_p}/2,
\end{equation}
where the first inequality follows from \eqref{i:fre2}.

Therefore  by Remark~\ref{r:sym}, we have that 
\begin{equation}\label{e:deco}
\tilde \Delta_{Y_p} = (\tilde \Delta_{Y_p}|_{\tilde{X}_{Y_p}\cap [0,K)} * \tilde \Delta_{Y_p}|_{\tilde{X}_{Y_p}\cap [K,N_{Y_p}/2]})^{\Sym}.
\end{equation}
Notice that $K \in S(r_-,r_+) \cap [0,2r_- + (p-3)r_+]$ and $K \neq (p-2)r_+$, so $K \in \tilde{S}_{Y_p}$ by Proposition~\ref{prop:reducedstructure}.  Therefore, $\tilde \Delta_{Y_p}|_{\tilde{X}_{Y_p}\cap [K,N_{Y_p}/2]}$ is an abstract delta sequence. We define
\begin{align} 
\Delta_{Z_p} &= \tilde{\Delta}_{Y_p}|_{\tilde{X}_{Y_p} \cap [0,K)} \\  
\Delta_{W_p} &= \tilde \Delta_{Y_p}|_{\tilde{X}_{Y_p}\cap [K,N_{Y_p}/2]}. \label{e:Wp} 
\end{align}
It is clear that $\Delta_{Z_p}$ is an abstract delta sequence, since $\tilde{\Delta}_{Y_p}$ is. 

\begin{lemma}\label{l:sink} 
The abstract delta sequence $\Delta_{Z_p}$ is sinking.  
\end{lemma}

\begin{proof}
We must verify all three properties in Definition~\ref{d:sink} for $\Delta_{Z_p}$.  Recall that $\tilde \Delta_{Y_p}$ is in reduced form, and consequently so is the restriction $\Delta_{Z_p}$.  Therefore, by definition of $\Delta_{Z_p}$, the last element of the delta sequence $\Delta_{Z_p}$ must have a negative value or else $\tilde \Delta_{Y_p}$ would contain two positive values in a row (namely the last element of $\tilde{X}_p \cap [0,K)$ and $K$), contradicting $\tilde \Delta_{Y_p}$ being in reduced form.   This establishes  Definition~\ref{d:sink}\eqref{i:sink1}, for $\Delta_{Z_p}$.

Before proceeding further, we set up notation. Throughout this proof, we will denote predecessors and successors taken with respect to $\tilde{X}_{Y_p}$ by a tilde decoration, and those taken with respect to $X_{Y_p}$ will not receive a tilde decoration.  Note that by the discussion in Section \ref{subsec:success}, we have
\begin{equation}\label{e:succesred}
\suc_+(x) \leq \widetilde{\suc}_+(x)\text{ for every }x \in \tilde{X}_{Y_p}.
\end{equation}

Next we show 
\begin{equation}\label{e:sink1}
\tilde \Delta_{Y_p}(x) \leq - \tilde \Delta_{Y_p}(\widetilde{\suc}_-(x)) \text{ for all  }x \in \tilde{S}_{Y_p} \cap [0,K),
\end{equation}
which will prove that the monotonicity condition in Definition~\ref{d:sink}\eqref{i:sink2}  holds for $\Delta_{Z_p}$. Let $x\in\tilde{S}_{Y_p} \cap [0,K)$. Then $x\in S(r_-,r_+)\cap [0,(p-3)r_+]$ by \eqref{e:Kest} and Proposition~\ref{prop:reducedstructure}\eqref{i1:reducedstructure}. Writing $x=ar_-+br_+$, we see that  $\tilde \Delta_{Y_p}(x) = \min\{a,b\} + 1$ by Proposition~\ref{prop:reducedstructure}\eqref{i2:reducedstructure}.  Let $y = (p-a-1)r_- + (p-b-3) r_+$.  By Lemma~\ref{lem:intervals}\eqref{shufflei:1}, and \eqref{e:succesred}, we have
\[
x < N_{Y_p} - y < \suc_+(x) \leq \widetilde{\suc}_+(x).    
\]
Note that since $x\in S(r_-,r_+)\cap [0,(p-3)r_+]$, we have that $a + b \leq p-3$.  Therefore, $(p-a-1)\geq 0$ and $(p-b-3) \geq 0$, and thus $N_{Y_p} - y \in Q_{Y_p}$.  By Proposition~\ref{prop:reducedstructure}\eqref{i3:reducedstructure}, $$\tilde \Delta_{Y_p}(\widetilde{\suc}_-(x)) \leq - \min\{p-a-1,p-b-3\} -1.$$ 
Therefore, to prove \eqref{e:sink1}, it suffices to show that 
\begin{equation}\label{eq:ab}
\min\{a,b\} \leq \min\{p-a-1,p-b-3\}. 
\end{equation}
However, we have seen that $a + b \leq p-3$.  Hence, $a \leq p-b-3$ and $b \leq p-a-3$, proving \eqref{eq:ab}.   

It remains to show that $\Delta_{Z_p}$ satisfies  Definition~\ref{d:sink}\eqref{i:sink3}. We observe that the last positive value of $\Delta_{Z_p}$ occurs at $\widetilde{\pre}_+(K) = \xi r_- + (\xi-1)r_+$ by Lemma~\ref{lem:sem2} and Proposition~\ref{prop:reducedstructure}\eqref{i1:reducedstructure}.  Thus, $\widetilde{\suc}_-(\xi r_- + (\xi-1)r_+)$ is the maximal element of $Z_p$. Hence we must prove 
\begin{equation}\label{e:lastZ}
\tilde{\Delta}_{Y_p}(\xi r_- + (\xi-1)r_+) < -\tilde{\Delta}_{Y_p}(\widetilde{\suc}_-(\xi r_- + (\xi-1)r_+)).
\end{equation}
We have $\tilde{\Delta}_{Y_p}(\xi r_- + (\xi-1)r_+) = \xi$ by Proposition~\ref{prop:reducedstructure}\eqref{i2:reducedstructure}. 
On the other hand, by Lemma~\ref{lem:intervals}\eqref{shufflei:1},   
\[
\xi r_- + (\xi-1)r_+ < N_{Y_p} - (p-\xi-1)r_- - (p-\xi-2)r_+  < \suc_+(\xi r_- + (\xi-1)r_+) \leq K.  
\]
Hence we also have 
\[
-\tilde \Delta_{Y_p}(\widetilde{\suc}_-(\xi r_- + (\xi-1)r_+)) \geq \min\{p-\xi-1,p-\xi-2\} + 1 = p-\xi-1
\]
by Proposition~\ref{prop:reducedstructure}\eqref{i3:reducedstructure}.  Therefore in order to show \eqref{e:lastZ}, it suffices to prove that $p-\xi-1>\xi$.  This is clear since $p = 2\xi+2$.
\end{proof}

\begin{lemma}\label{l:credec}
As abstract delta sequences, $\Delta_{W_p} \cong \Delta_{C_p}$, where $\Delta_{C_p}$ is the abstract delta sequence from Definition~\ref{d:cre}, and $\Delta_{W_p} $ is defined as in \eqref{e:Wp}.  
\end{lemma}
\begin{proof}
We would like to see that $\Delta_{W_p}$ agrees with $\Delta_{C_p}$ as abstract delta sequences.  To do this, we explicitly compute  $\Delta_{W_p}$. 

We first list all the elements of $\tilde{S}_{Y_p}\cap[K,N_{Y_p}/2]$. By Lemma~\ref{lem:sem2}, \eqref{i:fre2}, and \eqref{e:Kest}, it follows that $K<N_{Y_p}/2<(p-2)r_+$.  Now by Proposition~\ref{prop:reducedstructure}\eqref{i1:reducedstructure} we have $\tilde{S}_{Y_p}\cap[K,N_{Y_p}/2]=S(r_-,r_+)\cap[K,N_{Y_p}/2]$. Then using Lemma~\ref{lem:sem2}, we see that
\begin{align}
\nonumber \tilde{S}_{Y_p}\cap[K,N_{Y_p}/2]=& \{ (\xi-1)r_-+\xi r_+,\; (\xi-2)r_-+(\xi+1)r_+,\dots, \\
\label{e:SY}&\;r_-+(2\xi-2)r_+, (2\xi-1)r_+, 2\xi r_-,(2\xi-1)r_-+r_+,\dots, \\
\nonumber  & \qquad \qquad\; (\xi+2)r_-+(\xi-2)r_+, \; (\xi+1)r_-+(\xi-1)r_+ \} .
\end{align}
\noindent In order to verify that the last element in the above sequence is as indicated, we must show
\begin{align}
\label{i:xiNK<} (\xi+1)r_-+(\xi-1)r_+&<N_{Y_p}/2, \text{ and}\\
\label{i:xiNK>} \xi r_-+\xi r_+&>N_{Y_p}/2.
\end{align}
These two inequalities follow from \eqref{i:fre1} and \eqref{i:fre2} respectively, since $2\xi+2 = p$.  (For the first inequality, we use that $pr_- + (p-4)r_+ < (p-1)r_- + (p-3)r_+$.)  Hence \eqref{e:SY} holds.  

Similarly, we determine the elements of $\tilde{S}_{Y_p}\cap[N_{Y_p}/2, N_{Y_p}-K]$ so we may determine $\tilde{Q}_{Y_p}\cap[K,N_{Y_p}/2]$. Note that by Lemma~\ref{lem:sem2} and \eqref{i:fre2}, we have $N_{Y_p}/2<N_{Y_p}-K<2r_-+(p-3)r_+$ for $p \geq 4$.  By Proposition~\ref{prop:reducedstructure}\eqref{i1:reducedstructure}, we have $\tilde{S}_{Y_p}\cap[N_{Y_p}/2, N_{Y_p}-K]=S(r_-,r_+)\cap[N_{Y_p}/2, N_{Y_p}-K]\setminus \{2\xi r_+ \}$. Then using Lemma~\ref{lem:sem2}, \eqref{i:xiNK<}, and \eqref{i:xiNK>}, we see that
\begin{align}
\label{e:SY2} \tilde{S}_{Y_p}\cap[N_{Y_p}/2, N_{Y_p}-K]=& \{ \xi r_-+\xi r_+,\; (\xi-1)r_-+(\xi+1)r_+,\dots,\;r_-+(2\xi-1)r_+,  \\
\nonumber &(2\xi+1)r_-,2\xi r_-+r_+,\dots,\; (\xi+3)r_-+(\xi-2)r_+, \; (\xi+2)r_-+(\xi-1)r_+ \}.
\end{align}
\noindent Observe that we have purposely omitted $2\xi r_+$.  To see that the last element in the above sequence is as written, we must show
\begin{align}
(\xi+2)r_-+(\xi-1)r_+&<N_{Y_p}-K, \text{ and}\\
(\xi+1)r_-+\xi r_+&>N_{Y_p}-K.
\end{align}
Again these  inequalities follow from \eqref{i:fre1} and \eqref{i:fre2} respectively. Hence \eqref{e:SY2} holds.

By the  model for the reduced form of delta sequences described in Section~\ref{subsec:success}, the definition of $\Delta_{Y_p}$, and Proposition~\ref{prop:reducedstructure}, we have 
\begin{align}
\label{e:QY} \tilde{Q}_{Y_p} \cap [K,N_{Y_p}/2] &=\{N_{Y_p}-x \mid x\in \tilde{S}_{Y_p}\cap[N_{Y_p}/2, N_{Y_p}-K]\} \\
\nonumber &= \{N_{Y_p} - x \mid x \in S(r_-,r_+) \cap [N_{Y_p}/2,N_{Y_p}-K], x \neq 2\xi r_+\}.
\end{align}

We need to detect the positions of the elements of  $\tilde{Q}_{Y_p} \cap [K,N_{Y_p}/2]$ relative to the elements of $ \tilde{S}_{Y_p}\cap[K,N_{Y_p}/2]$. We shall employ the following inequalities
\begin{align}
\label{e:order1}  (\xi-1-j)r_-+(\xi+j)r_+ <N_{Y_p}-(\xi+2+j)r_--(\xi-1-j)r_+,\; j=0,\dots,\xi-1,\\
\label{e:order2}  N_{Y_p}-(\xi+2+j)r_--(\xi-1-j)r_+ <(\xi-2-j)r_-+(\xi+1+j)r_+,\; j=0,\dots,\xi-2,\\
\label{e:order3}  N_{Y_p}-jr_--(2\xi-j)r_+ <(2\xi-j)r_-+jr_+,\; j=0,\dots,\xi-1,\\
\label{e:order4}  (2\xi-j)r_-+jr_+ < N_{Y_p}-(j+1)r_--(2\xi-j-1)r_+,\; j=0,\dots,\xi-1.
\end{align}
These inequalities again follow from \eqref{i:fre1} and \eqref{i:fre2}. Finally, by Lemma \ref{lem:sem2}, \eqref{i:fre1}, and \eqref{i:fre2}, we point out that
\begin{equation}\label{e:order5}
(2\xi-1)r_+<N_{Y_p}-(2\xi+1)r_-<N_{Y_p}-2\xi r_+<2\xi r_-.
\end{equation}
Now it follows from \eqref{e:SY}, \eqref{e:SY2}, \eqref{e:QY}, \eqref{e:order1}, \eqref{e:order2}, \eqref{e:order3}, \eqref{e:order4}, and \eqref{e:order5} that the sequence $\tilde X_{Y_p} \cap [K,N_{Y_p}/2]$ is given by:
\begin{align}
\label{e:XY}
 \tilde X_{Y_p} \cap [K,N_{Y_p}/2]=& \{ (\xi-1)r_-+\xi r_+,\; N_{Y_p}-(\xi+2)r_--(\xi-1)r_+,\;\\
 \nonumber &(\xi-2)r_-+(\xi+1)r_+,\; N_{Y_p}-(\xi+3)r_--(\xi-2)r_+,\dots,\\
 \nonumber &\;r_-+(2\xi-2)r_+, \; N_{Y_p}-2\xi r_--r_+,(2\xi-1)r_+, \\
\nonumber & N_{Y_p}-(2\xi+1)r_-,\; 2\xi r_-, \; N_{Y_p}-r_--(2\xi-1)r_+,\; (2\xi-1)r_- + r_+,\dots,\\
\nonumber &\; N_{Y_p}-(\xi-1)r_--(\xi+1)r_+,\; (\xi+1)r_-+(\xi-1)r_+,\; N_{Y_p}-\xi r_--\xi r_+ \} . 
\end{align}
\noindent Again, note that $N_{Y_p} - 2\xi r_+$ is deliberately excluded from the above list, since it is not an element of $\tilde{X}_{Y_p}$.  

It remains to see that the values of $\tilde{\Delta}_{Y_p}$ on the above sequence are the same as the values of $\Delta_{C_p}$. By Proposition~\ref{prop:reducedstructure}\eqref{i2:reducedstructure}, since $N_{Y_p}/2 < (p-2)r_+$, we have that 
\begin{equation}\label{e:deltaonSY}
\tilde{\Delta}_{Y_p}(cr_-+dr_+)=\min \{ c,d\} +1,\text{ for } cr_-+dr_+\in \tilde S_{Y_p} \cap [K,N_{Y_p}/2].
\end{equation} 
\noindent Moreover for every $N_{Y_p}-cr_--dr_+\in \tilde Q_{Y_p} \cap [K,N_{Y_p}/2]$ such that $cr_-+dr_+\neq (2\xi+1)r_-$, we have 
\begin{equation}\label{e:deltaonQY}
\tilde{\Delta}_{Y_p}(N_{Y_p}- cr_--dr_+)=-\tilde{\Delta}_{Y_p}(cr_- + dr_+) = -\min \{ c,d\} -1,
\end{equation} 
by Proposition~\ref{prop:reducedstructure}, since as we observed earlier, $cr_- + dr_+ < N_{Y_p} - K < 2r_- + (p-3)r_+$.
Also, by Proposition~\ref{prop:reducedstructure}, we have $N_{Y_p} - (2\xi+1)r_- \in \tilde{Q}_{Y_p}$, and further,   
\begin{equation}\label{e:deltaonSp}
\tilde{\Delta}_{Y_p}(N_{Y_p}-(2\xi+1)r_-)= - \tilde{\Delta}_{Y_p}((2\xi+1)r_-) = -2.
\end{equation} 
Computing the value of $\tilde{\Delta}_{Y_p}$ at each element of the sequence \eqref{e:XY} using  \eqref{e:deltaonSY}, \eqref{e:deltaonQY}, and \eqref{e:deltaonSp}, and comparing the result with Definition~\ref{d:cre}, we  see that $\Delta_{W_p}$ agrees with $\Delta_{C_p}$. This is what we wanted to show.
\end{proof}

Having collected all the necessary ingredients, the proof of Lemma~\ref{l:dec} now follows from  \eqref{e:deco}, Lemma~\ref{l:sink}, and Lemma~\ref{l:credec}.  As discussed, this was the remaining piece needed to prove Proposition~\ref{prop:rank}, and consequently Theorem~\ref{thm:computation}.   

%% file: universalspheres.tex
\section{Knot surgeries in other manifolds}\label{sec:extra}
In the proof of Theorem~\ref{thm:obstruction}, the only thing special about $S^3$ is that it is an integer homology sphere L-space (i.e., $\HFred = 0$) and that $d(S^3) = 0$.  The following theorem is a slight generalization.

\begin{theorem}\label{thm:generalobstruction}
Let $Y$ and $Y'$ be oriented integer homology spheres.  Suppose that $\HFred(Y) = 0$ and $d(Y') \leq d(Y) - 8$.  Then, if $Y'$ is obtained by surgery on a knot in $Y$, then there is a non-trivial element of $\HFred(Y')$ in degree $d(Y)$ which is not in the kernel of $U$.  
\end{theorem}
\begin{proof}
The proof is the same as Theorem~\ref{thm:obstruction}, where one only difference is that one has to incorporate the $d$-invariant of $Y$ into some of the statements.  The main observation is that for $n > 0$, we have the more general formula $d(Y_{1/n}(K)) = d(Y) - 2V_0$, where $V_0$ is defined analogously for $Y$ and $K$ as for a knot in $S^3$.  This follows by repeating the arguments in \cite[Proposition 1.6]{NW} for a knot in an integer homology sphere L-space.        
\end{proof}

\begin{proof}[Proof of Theorem~\ref{thm:extra}]
Orient $\Sigma(2,3,5)$ such that it is the boundary of a negative-definite plumbing.  In this case, $d(\Sigma(2,3,5)) = 2$.  Let $Z$ be an integer homology sphere.  We will show that for $k \gg 0$, the manifold $\#_k \Sigma(2,3,5)$ is not surgery on a knot in $Z$, regardless of orientation of $Z$.  

Fix an orientation on $Z$.  Recall that $\HFred(Z)$ is finite-dimensional over $\F$.  Therefore, we may define an integer $n_Z$ by 
\[
n_Z = | \max \{ s \mid \HFred_s(Z) \neq 0\}|. 
\]
Choose $k>0$ such that 
\[
2k \geq \max\{d(+Z),d(-Z)\} + \max\{n_{+Z},n_{-Z}\} + 8. 
\]
Due to the additivity of $d$ under connected-sums \cite[Theorem 4.3]{OSGraded}, we have 
\[
d(\#_k \Sigma(2,3,5)) = 2k \geq d(\pm Z) + n_{\pm Z} + 8.
\]
Therefore, $d(\pm Z) \leq d(\#_k \Sigma(2,3,5)) - 8$.  Furthermore, by construction there is no element of $\HFred(\pm Z)$ in degree $d(\#_k \Sigma(2,3,5))$.  Therefore, by Theorem~\ref{thm:generalobstruction}, neither $Z$ nor $-Z$ can be expressed as surgery on a knot in $\#_k \Sigma(2,3,5)$.  Consequently, $\#_k \Sigma(2,3,5)$ cannot be surgery on a knot in $Z$, regardless of orientation.                  
\end{proof}

We conclude by pointing out that Theorem~\ref{thm:generalobstruction} can also be extended to statements about $p/q$-surgery where $|p| \geq 2$.  One can then apply the same arguments as in Theorems~\ref{thm:obstruction} and ~\ref{thm:extra} to show that if $Z$ has cyclic first homology, $\#_k \Sigma(2,3,5)$ is not surgery on a knot in $Z$ for $k$ large.  Finally, the analogous statement when $Z$ has non-cyclic homology is trivial.  Thus, in conclusion, for any three-manifold $Z$, there exist infinitely many integer homology spheres which are not surgery on a knot in $Z$.  
 
       
